\documentclass[11pt,english]{article} 
\usepackage[utf8]{inputenc}
\usepackage[T1]{fontenc}
\usepackage{lmodern}
\usepackage[a4paper]{geometry}
\usepackage{babel}
\usepackage{amsmath}
\usepackage{cleveref}

\usepackage{nag}
\usepackage{xcolor}
\usepackage{epsfig,rotating}
\usepackage{amsmath}
\usepackage{mathtools}
\usepackage{amsthm}
\usepackage{amssymb}
\usepackage{graphicx}
\usepackage{caption}
\usepackage{enumitem}

                \newcommand{\R}{\mathbb R}

\newcommand{\D}{\partial}               
\newcommand{\s}{\sigma}                 \newcommand{\g}{\gamma} 
\newcommand{\e}{\varepsilon}            \newcommand{\G}{\Gamma} 
                 \def\l{\lambda}
\newcommand{\sph}{\mathbb S}                        
\newcommand{\comp}{{\small{\circ}}}     \def\a{\alpha}
               \newcommand{\E}{{\mathcal E}}

\newcommand{\N}{{\mathbb N}}            \newcommand{\RP}{{\mathbb R\mathrm P}}

\def\t{\tau}
\def\hat{\widehat}                        
                 
\newcommand{\Abf}{\boldsymbol{A}}

\newcommand{\Qbf}{\boldsymbol{Q}}               \newcommand{\lbf}{\boldsymbol{\lambda}} 
\newcommand{\Qp}{\boldsymbol{Q}_{\scriptscriptstyle{\mathrm{P}}}}  \newcommand{\Qs}{\boldsymbol{Q}_{\scriptscriptstyle{\mathrm{S}}}}
               \newcommand{\nbf}{\boldsymbol{n}}
\newcommand{\ubf}{\boldsymbol{u}}

\begin{document}

\captionsetup[figure]{labelfont={bf},labelformat={default},labelsep=period,name={\small Fig.}}


\numberwithin{equation}{section}                
\theoremstyle{plain}
\newtheorem{theorem}{\bf Theorem}
\newtheorem*{theorem*}{\bf Theorem}
\newtheorem{teorema}{\bf Theorem}
\def\theteorema{\Roman{teorema}}
\newtheorem{lemma}{\bf Lemma}
\newtheorem*{lemma*}{\bf Lemma}
\newtheorem*{mlemma}{\bf Main Lemma}

\newtheorem{proposition}{\bf Proposition}
\newtheorem*{proposition*}{\bf Proposition}
\newtheorem{corollary}{\bf Corollary}
\newtheorem*{corollary*}{\bf Corollary}
\newtheorem*{conjecture}{\bf Conjecture}
\newtheorem*{fact*}{\bf Fact}
\newtheorem{property}{\bf Property}


\theoremstyle{definition}
\newtheorem{definition}{\bf Definition}
\newtheorem*{definition*}{\bf Definition}

\newtheorem{example}{\bf Example}
\newtheorem*{example*}{\bf Example}
\newtheorem{Example}{\bf Example}[section]
\theoremstyle{remark}
\newtheorem*{remark*}{\bf Remark}
\newtheorem{remark}{\bf Remark}
\newtheorem*{problem}{\bf Problem}
\newtheorem*{exercise*}{\bf Exercise}


\title{\vspace*{-2cm} On the Geometry of the Quadratic Maps}

\author{By Ricardo\,Uribe-Vargas\footnote{IMB-UMR5584 CNRS, Universit\'e Bourgogne Franche-Comt\'e and 
  Laboratory Solomon Lefschetz UMI2001 CNRS, 
  Universidad Nacional Autonoma de M\'exico. \textit{r.uribe-vargas@u-bourgogne.fr}}}
  

\date\empty                     
\maketitle



\let\h\theta
\let\t\tau

\begin{abstract}
\noindent
To analyse a homogeneous quadratic map $\Qbf:\R^k\to\R^\ell$, the critical value variety $\Sigma_{\Qbf}$ of its restriction to the unit sphere, 
$\Qbf_{|\sph^{k-1}}$, is very useful but difficult to handle. 
We prove its dual hypersurface $\mathcal{F}_{\Qbf}$ 
in $\R^\ell\sqcup\RP_\infty^{\ell-1}$ is of degree $k$ (providing its equation) and is \textit{hyperbolic}: 
every line though the origin in $\R^\ell$ cuts $\mathcal{F}_{\Qbf}$ at $k$ real points, 
counting multiplicities. These facts imply: 
\textbf{1}. $\mathcal{F}_{\Qbf}$ is the union of $k/2$ nested  topological hypersheres if $k$ is even, or the union of $(k-1)/2$ nested  topological hypersheres and one topological projective hyperplane if $k$ is odd.   
\textbf{2}.\,The Toeplitz-Haussdorf Theorem on the range of $\Qbf_{|\sph^{k-1}}$ if $k>\ell$. 
\textbf{3}. The dual version of the hyperbolicity of  $\mathcal{F}_{\Qbf}$ is: the critical value variety $\Sigma_{\Qbf}$ is of class $k$. 
\textbf{4}. For generic $\Qbf:\R^k\to\R^k$, $\Sigma_{\Qbf}$ contains $\lfloor\frac{k-1}{2}\rfloor$ collections of 
ellipsoids of codimension $2$ 
consisting of  degenerate 
parabolic points. 
\textbf{5}. If $k>\ell$, the parabolic set in $\Sigma_{\Qbf}$
(if any) is a union of hyperplane closed discs bounded by ellipsoids.
\textbf{6}. The above ellipsoids (items $4$ and $5$) 
are in bijection with the isolated singular points in $\mathcal{F}_{\Qbf}$. 
\end{abstract}
\smallskip
. 
%


\section*{Introduction}
Homogeneous quadratic maps appear in many areas of Mathematics and Physics, being the jets of degree $2$ 
of the maps $\R^k\to\R^\ell$ at its points of rank $0$. For example, for a submanifold $M^k$ in $\R^{k+\ell}$, given as the graph of a map $f:\R^k\to\R^\ell$, at a point $p\in M$,  
its differential invariants of order two at $p$ are determined by the $2$-jet of $f$ at $O_{\R^k}$ (cf. \cite{Little}). 
The geometry, topology and applications of homogeneous quadratic maps have been studied 
by Wall, Lopez de Medrano, Gitler and many others, cf. \cite{Wall,Santiago1,Santiago-Samuel} (see more references in \cite{Santiago-Samuel}). 
\smallskip 

We start with a homogeneous quadratic map $\Qbf:\R^k\to\R^\ell$, 
\[t\mapsto Q_1(t)\nbf_1+\cdots+Q_\ell(t)\nbf_\ell\,,\]
where $t=(t_1,\ldots,t_k)$ orthonormal coordinates in $\R^k$, $\nbf_1,\ldots,\nbf_\ell$ is an orthonormal basis of $\R^\ell$ and 
$Q_r:\R^k\to\R$ are quadratic forms.  \ 
It is written also as 
\begin{equation}\label{eq:homogeneous-Rk_Rn}
	\qquad \Qbf:t\mapsto\sum\,t_it_j\Abf_{ij}\,, \quad \mbox{ with } \Abf_{ij}=\Abf_{ji}\in\R^\ell\,,
\end{equation}
where $\,\Abf_{ij}:=a_{ij}^1\nbf_1+\cdots+a_{ij}^\ell\nbf_\ell$,\,
and $\left\{a^r_{ij}\right\}$ is the matrix of the quadratic form $Q_r$. 

\noindent
\textbf{\textit{Restriction of $\Qbf$ to $\sph^{k-1}$.}} 
The homogenity, $\Qbf(\l t)=\l^2\Qbf(t)$, implies that $\Qbf$ is determined by its restriction to the unit hypersphere $\sph^{k-1}$. We shall consider it:
\[\qquad \Qbf_{_\mathrm{S}}:=\Qbf_{|\sph^{k-1}}:\sph^{k-1}\to\R^\ell\,, \qquad 
\mbox{where }\, \sph^{k-1}=\{t\in\R^k:\|t\|=1\}\,.\]
Moreover, since $\Qs(-t)=\Qs(t)$, we also consider $\Qs$ as a map from $\RP^{k-1}$,   
\begin{equation}\label{eq:formula-quad-map}
 \Qbf_{_\mathrm{P}}:\RP^{k-1}\to\R^\ell\,, \qquad \Qp([t])=\frac{1}{\|t\|^2}\Qbf(t)\,,
\end{equation}
where $[t]\in\RP^{k-1}$ is the class of all non-zero multiples of $t\in\R^k\smallsetminus\{0\}$.   

\section{Front and Focal of a Quadratic map $\Qbf:\R^k\to\R^\ell$}
\subsection*{A \, Front (Critical Value Set)}
If $k-1>\ell$, the fibres $\Qp^{-1}(\omega)$ of $\Qp:\RP^{k-1}\to\R^\ell$ (or of $\Qs$) are intersections of quadrics, whose topological type 
(of great interest, cf. \cite{Santiago-libro} and its references) changes if the value $\omega$ crosses the ``critical value set''.
\medskip

\noindent
\textit{\textbf{{\small Critical Values}}}. 
A \textit{critical point} of a map $f:M^k\to N^\ell$ is a point $p\in M$ at which the differential $d_pf$ 
is not ``onto''\,: $\mathrm{rank} \ d_pf<\ell$. \\
A \textit{critical value} of $f$ is its value at a critical point. \\
(If $k<\ell$, all points of $M^k$ are critical and the critical value set is $f(M)\subset N$.)
\medskip

\noindent
\textit{\textbf{{\small Front}}}. We call \textit{front} of the quadratic map $\Qbf:\R^k\to\R^\ell$, the critical value set of the map $\Qp$ (or of $\Qs$), and denote it by $\Sigma_{\Qbf}$.  Notice that $\Qs$ covers twice $\Sigma_{\Qbf}$. 

{\footnotesize
\begin{remark*}
It turns out that the critical value set $\Sigma_{\Qbf}$ of $\Qp$ is the front of a Legendrian map and we shall use some front features of $\Sigma_{\Qbf}$ (\cite{avg}). Whence our terminology. 
\end{remark*}
}
\begin{example*}[Veronese Manifolds]
A \textit{Veronese $(k-1)$-manifold} is the image of a homogeneous quadratic map 
$\Qp:\RP^{k-1}\to\R^{\frac{1}{2}k(k+1)}$, where the vectors $\Abf_{ij}$ ($1\leq i\leq j\leq k$) are independent - we note it $V_{k-1}$. 
One easyly checks that 
\textit{$V_{k-1}$ is an inflection-free embedding of $\RP^{k-1}$ in 
a hyperperplane of $\R^{\frac{1}{2}k(k+1)}$, and
lies in a hyperellipsoid of $\R^{\frac{1}{2}k(k+1)}$}. Moreover, 
\textit{all Veronese $(k-1)$-manifolds are affinely equivalent in $\R^{\frac{1}{2}k(k+1)}$}. 

Thus, a Verenese $(k-1)$-manifold $\Qp(\RP^{k-1})$ is a front: $\Qp(\RP^{k-1})=\Sigma_{\Qbf}$. \\
For a homogeneous quadratic map $\Qbf:\R^k\to\R^\ell$ with 
$k-1<\ell<\frac{1}{2}k(k+1)$, the front $\Sigma_{\Qbf}=\Qp(\RP^{k-1})$ is the projection of a Verenese $(k-1)$-manifold.
\end{example*}

\begin{example*}[Steiner Roman Surface]
Take the quadratic map $\Qbf:\R^3\to\R^3$, given by 
\begin{equation}
\Qbf(t)=(2t_2t_3, 2t_1t_3, 2t_1t_2)\,.
\end{equation} 
Sice $k-1<\ell$, $\Sigma_{\Qbf}=\Qp(\RP^2)$ and the front is the 
\textit{Steiner Roman Surface} (Fig.\,\ref{fig:Steiner-Roman-Surface}): 
\begin{figure}[ht!]
	\centering
	\includegraphics[width=0.37\textwidth]{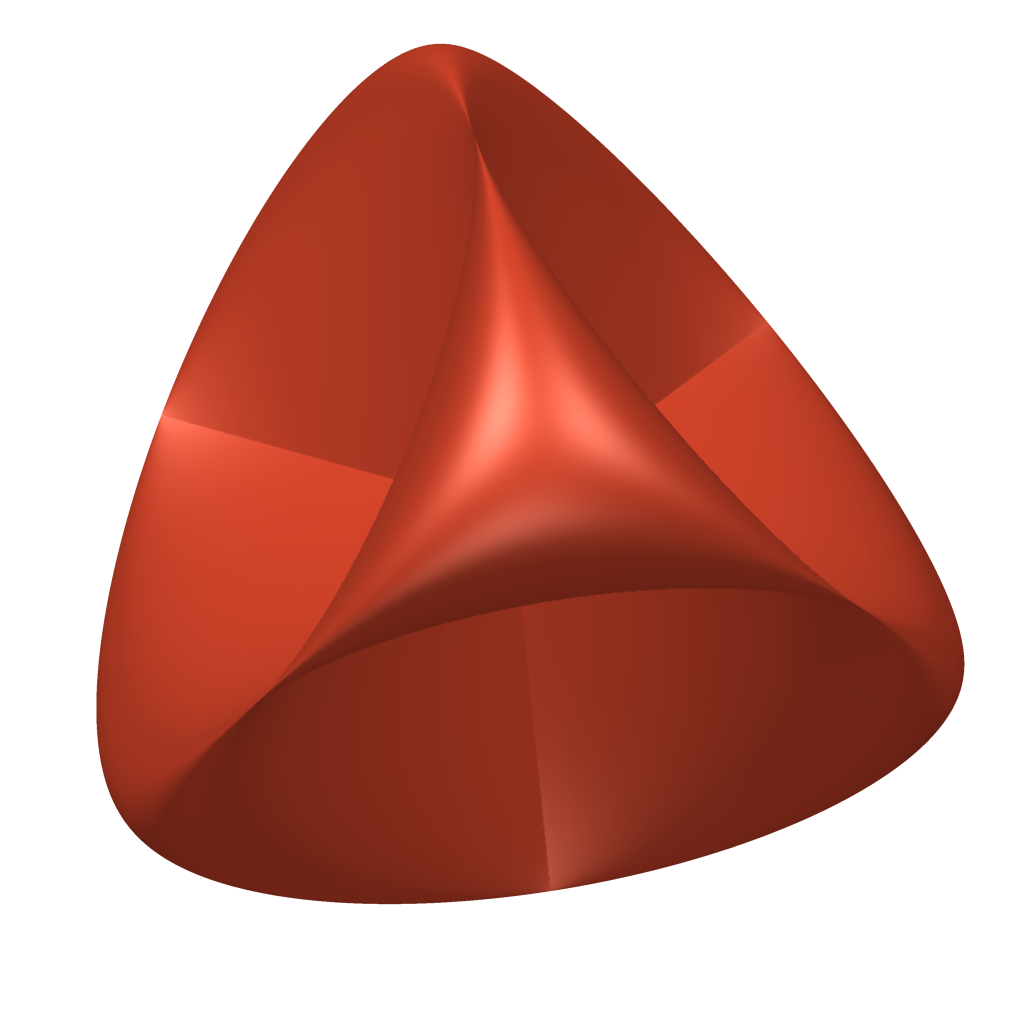}
	\caption{\small The Steiner Roman Surface is the front of the map $\Qbf(t)=(2t_2t_3, 2t_1t_3, 2t_1t_2)$.}
	\label{fig:Steiner-Roman-Surface}
\end{figure}
\end{example*}

\begin{example*}[Local Invariants of Submanifolds]
A $k$-dimensional submanifold $M^k$ in $\R^{k+\ell}$ is locally given in ``Monge form'' $f:\R^k\to\R^{k+\ell}$, $f(t)=(t;\frac{1}{2}\Qbf(t)+h.o.t.(t))$ at $p=f(0)$, 
where $\Qbf:\R^k\to\R^\ell$ is a homogeneous quadratic map.  

The image of $\Qp:\RP^{k-1}\to\R^\ell(=N_pM)$ is called \textit{curvature locus} of $M$ at $p=f(0)$. 
All 2nd order local invariants of $M$ at $p$ depend on the front $\Sigma_{\Qbf}$ which, in this context, is called \textit{curvature indicatrix}. 
The case of surfaces ($k=2$), in which $\Sigma_{\Qbf}$ is an ellipse, was wadely studied in the last $50$ years (see references in \cite{IRFRT}). 
\end{example*}

\begin{example*}[Diagonal Maps and Polytopes]
The \textit{diagonal maps} (fundamental to study intersections of quadrics \cite{Santiago-libro}) are those whose quadratic forms $Q_r$ 
contain only sums of squares $t_i^2$.
\textit{The front of a diagonal map is a polytope}. Example: if $\Qbf:\R^5\to\R^2$ with $\Qbf(t)=\sum_{i=1}^5t_i^2\Abf_{ii}$ 
(thus $\Qp:\RP^4\to\R^2$), then the front $\Sigma_{\Qbf}$ is a pentagon together with all its diagonal segments - all edges having multiplicity two. 
\end{example*}

\begin{figure}[ht!]
\centering
\includegraphics[width=1\textwidth]{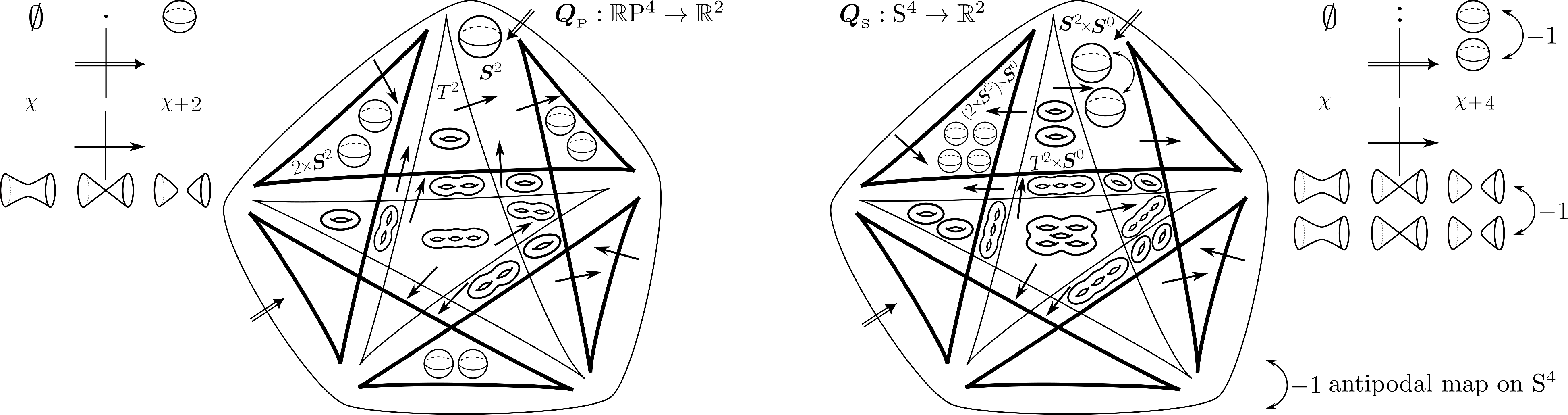} 
\caption{\small The fibres of $\Qp$ and $\Qs$ for the diagonal pentagonal map} 
\label{fig:pentagon-fibres}
\end{figure}

\noindent
\textbf{A Difficulty with the Front $\Sigma_{\Qbf}$}. If $k-1\geq \ell$, $\Sigma_{\Qbf}\neq\Qbf(\RP^{k-1})$ and, 
besides the case of diagonal maps, the front $\Sigma_{\Qbf}$ is not easy to compute and to handle. 
However, its dual hypersurface in $\RP^\ell$ (called ``focal'' below) has simple equation and nice properties. 
For example, for a generic small deformation of the above pentagon  diagonal map, it enabled us to guess the shape of the front $\Sigma_{\Qbf}$ and the fibres $\Qp^{-1}(\omega)$ and $\Qs^{-1}(\omega)$ (Fig.\,\ref{fig:pentagon-fibres}).  

In Fig.\,\ref{fig:pentagon-fibres}, the arrows indicate the local Morse transitions in which the Euler characteristic shifts by $+2$. Note that we can get the fibres of $\Qp$ from those of $\Qs$, invariant under the antipodal map in $\mathrm S^4$, by taking their quotient by the action of that map.

\subsection*{B \, Focal (Local Caustic) - a Degree $k$ Hypersurface} 
The maps $\Qbf$, $\Qp$ and $\Qs$ can be studied through 
the graph of $\frac{1}{2}\Qbf:\R^k\to\R^\ell$ - a paraboloid we denote by $\G_{\frac{1}{2}\Qbf}$ - 
near the origin $(O_{\R^k};O_{\R^\ell})$.    
Here we consider $\Qbf$ as the quadratic map associated to the ``Hessian'' 
of $\frac{1}{2}\Qbf$ at the origin, $\frac{\D^2}{\D t_i\D t_j}\left(\frac{1}{2}\Qbf\right)_{|_O}$. 

The centres of the hyperspheres in $\R^k\times\R^\ell$ which are tangent to $\G_{\frac{1}{2}\Qbf}$ at $(O_{\R^k};O_{\R^\ell})$ form the vectorial subspace 
$\left\{O_{\R^k}\right\}\times\R^\ell$, identified with $\R^\ell$. Its projective completion $\left\{O_{\R^k}\right\}\times(\R^\ell\sqcup\RP_{\infty}^{\ell-1})$, 
identified with $\R^\ell\sqcup\RP^{\ell-1}_\infty$, takes into account the tangent hyperspheres of infinite radius (the tangent hyperplanes). 
\medskip

\noindent
\textbf{\small\textit{Focal}}. The \textit{focal} 
of a homogeneous quadratic map $\Qbf:\R^k\to\R^\ell$, noted $\mathcal{F}_{\Qbf}$, 
is the projective completion, in $\R^\ell\sqcup\RP_{\infty}^{\ell-1}$, 
of the subset formed by the centres of the hyperspheres that have unusual tangency with the paraboloid $\G_{\frac{1}{2}\Qbf}$ at the origin. 

\begin{theorem}\label{th:Focal-Eq}
The focal $\mathcal{F}_{\Qbf}$ of a quadratic map $\Qbf:\R^k\to\R^\ell$ is the degree $k$ hypersurface in $\R^\ell\sqcup\RP_\infty^{\ell-1}$ given, in affine coordinates, 
by the equation 
\begin{equation}\label{eq:focal_set-equation}
\det
\begin{bmatrix}
	\langle \Abf_{11}\,, \, q\rangle -1 & \langle \Abf_{12}\,, \, q\rangle & \ldots & \langle \Abf_{1k}\,, \, q\rangle\\ 
	\langle \Abf_{21}\,, \, q\rangle & \langle \Abf_{22}\,, \, q\rangle -1 & \ldots & \langle \Abf_{2k}\,, \, q\rangle\\
	\vdots & \vdots & \ddots & \vdots \\
	\langle \Abf_{k1}\,, \, q\rangle & \langle \Abf_{k2}\,, \, q\rangle & \ldots & \langle \Abf_{kk}\,, \, q\rangle-1
\end{bmatrix} =0\,. 
\end{equation}
\end{theorem}

\begin{proof}
The centre of a hypersphere in $\R^k\times\R^\ell$ tangent to the paraboloid $\G_{\frac{1}{2}\Qbf}$ at $(O;O)$ has the form $(O;q)$. 
The equation of such hypersphere $S_q$ is
$\|(x;y)-(O;q)\|^2=\|(O;q)\|^2$, which is equivalent to the equation 
$\|x\|_{\R^k}^2+\|y\|_{\R^\ell}^2-2\langle y,q\rangle_{\R^\ell}=0$. 
So, writing 
\[
g_q(x;y):=\|x\|_{\R^k}^2+\|y\|_{\R^\ell}^2-2\langle y,q\rangle_{\R^\ell}\,,
\]
the hypersphere $S_q$, of equation $g_q(x;y)=0$, is osculating to  $\G_{\frac{1}{2}\Qbf}$ at $(O,O)$ iff the quadratic part (in $t$) 
	of the function \[g_q\comp\gamma(t)=\|t\|_{\R^k}^2+\|\tfrac{1}{2}\Qbf(t)\|_{\R^\ell}^2-
	2\langle \tfrac{1}{2}\Qbf(t),q\rangle_{\R^\ell}\] 
	is degenerate. That quadratic part, $\|t\|_{\R^k}^2-\langle \Qbf(t),q\rangle_{\R^\ell}$, is equal to 
	\[-\left(\sum_{i=1}^{k}t_i^2(\langle\Abf_{ii},q\rangle-1)+2\sum_{1\leq i<j\leq k}t_it_j\langle\Abf_{ij},q\rangle\right)\,.\]
	It degenerates iff the point $q\in\R^\ell\sqcup \RP^{\ell-1}_\infty$ satisfies  \eqref{eq:focal_set-equation}, in affine coordinates. 
\end{proof}

\begin{remark*}
The focal of $\Qbf$ can be seen as the projective completion of the set of points $q\in\R^\ell$ for which the quadratic form 
$Q_q(t):=\langle\,\Qbf(t)\,,\,q\,\rangle\,$ has $1$ as an eigenvalue. 
\end{remark*}

\subsection*{C \, Hyperbolicity of the Focal and its Global Configuration} 
\noindent
\textbf{\textit{\small Hyperbolicity}}. A degree $k$ variety $\mathcal{C}$ in $\RP^\ell$ is \textit{hyperbolic relative to a point} $p$ if it satisfies the 
``maximal intersection property'': 
\smallskip

\centerline{\sl Every line through $p$ cuts $\mathcal{C}$ at $k$ real points, counting multiplicities.}

\begin{theorem}\label{th:hyperbolicity-theorem}
Let $\Qbf:\R^k\to\R^\ell$ a homogeneaus quadratic map. 
	
\noindent
\textit{\textbf{a}}$)$ 
The focal hypersurface $\mathcal{F}_{\Qbf}$ is degree $k$-hyperbolic in $\R^\ell\sqcup\RP_\infty^{\ell-1}$ relative to 
$O$ in $\R^\ell$. 
\smallskip

\noindent
\textit{\textbf{b}}$)$
For any unit vector $\ubf\in\R^\ell$, the line generated by $\ubf$ 
meets $\mathcal{F}_{\Qbf}$ at the points $(1/\mu_\ell)\ubf$ $(\ell=1,\ldots,k)$, 
where each $\mu_\ell$ is an eigenvalue of the matrix
$(M_{\Qbf}^{\ubf})_{ij}:=\langle\,\Abf_{ij}\,,\,\ubf\,\rangle$. 
The eigenvectors $t\in\R^k$ associated to distinct eigenvalues are orthogonal. 
\end{theorem}

\begin{proof}
We write $P_{\Qbf}(q):=\det(M_{\Qbf}^q-\mathrm{Id}_k)$, where 
$(M_{\Qbf}^q)_{ij}=\langle\,\Abf_{ij}\,,\,q\,\rangle$. 
Thus equation \eqref{eq:focal_set-equation} is written as  $P_{\Qbf}(q)=0$. 
Let $\ubf$ be a unit vector in $\R^\ell$ and take the vectorial line generated 
by $\ubf$,  $L_{\ubf}=\{\l\ubf:\l\in\{\R\smallsetminus 0\}\}$. 
Its intersection points with the focal $\mathcal{F}_{\Qbf}$ are given by the solutions, in $\l$, of the equation 
$P_{\Qbf}(\l\ubf)=0$. Now, writing $\mu=1/\l$, the equation $P_{\Qbf}(\l\ubf)/\l^k=0$ is written as 
\begin{equation}\label{eq:polar-dual-egenvalues}
	\det
	\begin{bmatrix}
		\langle \Abf_{11}\,, \, \ubf\rangle -\mu & \langle \Abf_{12}\,, \, \ubf\rangle & \ldots & \langle \Abf_{1k}\,, \, \ubf\rangle\\ 
		\langle \Abf_{21}\,, \, \ubf\rangle & \langle \Abf_{22}\,, \, \ubf\rangle -\mu & \ldots & \langle \Abf_{2k}\,, \, \ubf\rangle\\
		\vdots & \vdots & \ddots & \vdots \\
		\langle \Abf_{k1}\,, \, \ubf\rangle & \langle \Abf_{k2}\,, \, \ubf\rangle & \ldots & \langle \Abf_{kk}\,, \, \ubf\rangle-\mu
	\end{bmatrix} =0. 
\end{equation}
The numbers $\mu$ satisfying \eqref{eq:polar-dual-egenvalues} are the eigenvalues of the matrix 
$(M_{\Qbf}^{\ubf})_{ij}:=\langle\Abf_{ij},\ubf\rangle$, all of which are real because 
$(M_{\Qbf}^{\ubf})_{ij}$ is symmetric. Thus, the line $L_{\ubf}$ intersects the focal 
$\mathcal{F}_{\Qbf}$ at $k$ real points, counting multiplicities. 
Distinct eigenvalues correspond to orthogonal unit eigenvectors $t\in\sph^{k-1}$. An eigenvalue $\mu$ 
corresponds to a point at $\RP^{n-1}_\infty$ iff $\mu=0$.  
\end{proof}

\begin{figure}[ht!]
\centering
\includegraphics [scale=0.18]{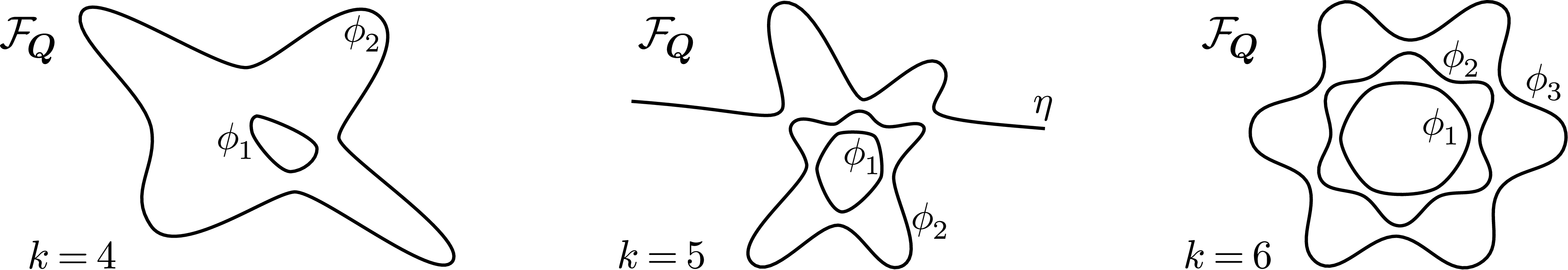}
\caption{\small Examples of focals 
$\mathcal{F}_{\Qbf}\subset\R^2\sqcup\RP^1_\infty$ of quadratic maps $\Qbf:\R^k\to\R^2$ for $k=4,5,6$.}
\label{fig:nested-smooth-spheres}
\end{figure}
\begin{figure}[ht!]
\centering
\includegraphics [scale=0.2]{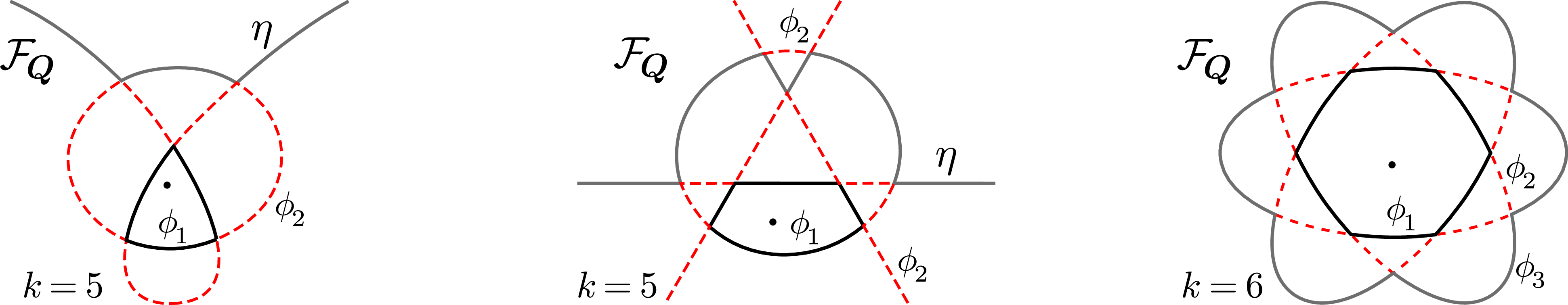}
\caption{\small Non smooth reducible focals $\mathcal{F}_{\Qbf}$. If $k=5$, $\mathcal{F}_{\Qbf}$ consists of two topological circles 
$\phi_1$, $\phi_2$ and one topological line $\eta$; if $k=6$, $\mathcal{F}_{\Qbf}$ consists of three topological circles $\phi_1$, $\phi_2$, $\phi_3$.}
\label{fig:nested-top-hyperspheres}
\end{figure}

The hyperbolicity of the focal $\mathcal{F}_{\Qbf}$ determines its topological configuration:

\begin{theorem}\label{th:topological-configurations}
	Let $\Qbf:\R^k\to\R^\ell$ a homogeneous quadratic map. 
	\smallskip
	
	\noindent
	\textit{\textbf{a}}$)$ The configuration of $\mathcal{F}_{\Qbf}$ in $\RP^\ell\,(\R^\ell\sqcup\RP^{\ell-1}_\infty)$ is as follows 
	{\rm (Figs.\,\ref{fig:nested-smooth-spheres} and \ref{fig:nested-top-hyperspheres})}$:$
	\smallskip
	
	\noindent 
	even $k:$ \, $\frac{k}{2}$ nested topological hyperspheres $\phi_1,\ldots,\phi_{k/2}$; 
	\smallskip
	
	\noindent
	odd $k:$  $\frac{(k-1)}{2}$ nested topological hyperspheres $\phi_1,\ldots,\phi_{(k-1)/2}$ and one topological hyperplane $\eta$. 
	\smallskip
	
\noindent
\textit{\textbf{b}}$)$ 
The domain bounded by the innermost topological hypersphere $\phi_1$ is convex. 
Moreover, relative to each point in that domain, $\mathcal{F}_{\Qbf}$ is hyperbolic and each $\phi_i$ is star-shaped. 
\end{theorem}
\begin{proof}
	\textit{\textbf{a}}) \underline{($k$ even)}: By the hyperbolicity of $\mathcal{F}_{\Qbf}$, each projective line $\lbf$ trough the origin $O$ of $\R^\ell$ 
	in $\R^\ell\sqcup\RP^{\ell-1}_\infty$ intersects $\mathcal{F}_{\Qbf}$ at $k$ real points counting multiplicities. 
	These $k$ points determine $k/2$ nested segments that contain $O$: 
	$I_1^{\lbf}\subseteq I_2^{\lbf}\subseteq\ldots\subseteq I_{k/2}^{\lbf}\,. $
	Taking all lines $\lbf$ through the origin $O$, the end-points of the segments $I_i^{\lbf}$ form a topological hypersphere $\phi_i$ ($i=1,\ldots,k/2$). 
	
	\noindent
	\underline{($k$ odd)}: Similarly, in each projective line $\lbf$ through $O$ there are $(k-1)/2$ nested segments that contain $O$, 
	$I_1^{\lbf}\subseteq I_2^{\lbf}\subseteq\ldots\subseteq I_{(k-1)/2}^{\lbf}\,,$
	giving rise to $(k-1)/2$ nested topological hyperspheres $\phi_1,\ldots,\phi_{(k-1)/2}$. 
	In each line $\lbf$, it remains one point $\tau_{\lbf}$ not contained in the interior of $I_{(k-1)/2}^{\lbf}$, but which may coincide with an end-point of $I_{(k-1)/2}^{\lbf}$. 
	Taking all lines $\lbf$ through $O$, the points $\tau_{\lbf}$ form a topological projective hyperplane $\eta$. 
	
	\smallskip
	
\noindent
\textit{\textbf{b}}) The configuration of the degree $k$ hypersurface $\mathcal{F}_{\Qbf}$, described in \textit{\textbf{a}},
imply that $\mathcal{F}_{\Qbf}$ is hyperbolic relative to each point in the domain bounded by the innermost topological hypersphere $\phi_1$. 
Thus, each line through a point in that domain cuts $\phi_1$ at exactly two points. This implies the convexity of that domain and that each $\phi_i$ is star-shaped relative to any point in that domain. 
\end{proof}

\noindent
\textbf{Note}. If a line through $O$ in $\R^\ell$ intersects $\mathcal{F}_{\Qbf}$ with multiplicity $j$ at a point, we shall consider that point as ``$j$ distinct 
points located at the same place''. 
With this convention, ``the nested hyperspheres $\phi_1,\ldots,\phi_{k/2}$ are disjoint, but may overlap at some (singular) points'' - Fig.\,\ref{fig:nested-top-hyperspheres}.
\medskip

\subsection*{D \, Hyperbolic singularities} 

The hyperbolicity of $\mathcal{F}_{\Qbf}$ restricts its possible singularities (characterised in \cite{Vainshtein}): 
\medskip

\noindent
\textbf{\textit{\small Local Hyperbolicity}} (Garding \cite{Garding}). A germ of hypersurface $f=0$ at $q\in\R^\ell$ is said to be \textit{hyperbolic} 
relative to a given line $\ell$ through $q$ if in a neighbourhood of $q$ the germ has a constant real intersection multiplicity with any real line sifficiently close to $\ell$. 
This local hyperbolicity was considered by Atiyah, Bott and Garding (see \cite{Garding,AtBG})

\begin{proposition}[\cite{Vainshtein}]
Any singularity of a hyperbolic hypersurface is locally hyperbolic. 
\end{proposition}

\begin{theorem}[\cite{Vainshtein}]\label{th:hyperbolic-singularities}
The stable locally hyperbolic singularities of hypersurfaces of dimension greater  
than $1$ are isolated and have multiplicity $2$, so that their equation is written in the form $f=0$, where \,$f\equiv x_1^2-\phi(x_2,\ldots,x_\ell)$, \,$\phi\geq 0$, \,$\phi(0)=0$. 
\end{theorem}

\begin{figure}[ht!]
	\centering
	\includegraphics [scale=0.2]{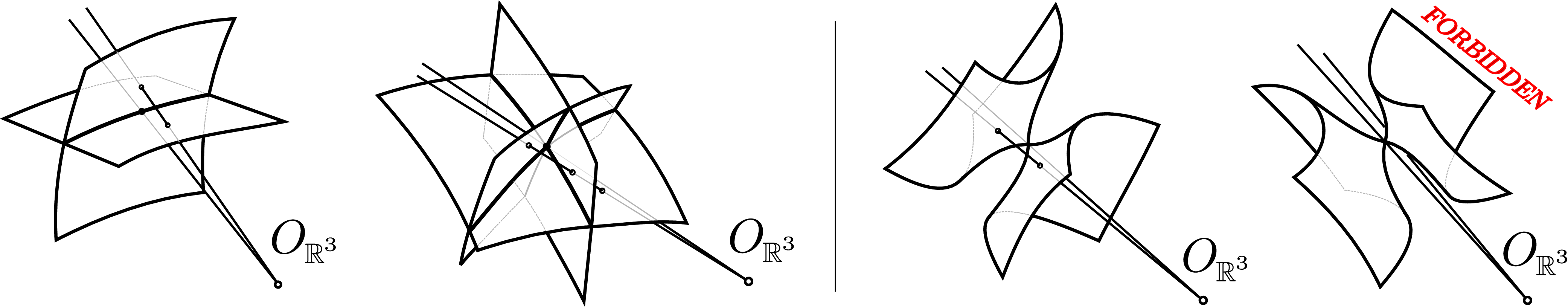}
	\caption{\small Left: some non isolated singularities of $\mathcal{F}_{\Qbf}$. \hspace*{1cm} Right: Morse singularities. }
	\label{fig:admissible-singul}
\end{figure}

The singularities of the focal $\mathcal{F}_{\Qbf}$ are hyperbolic relative to a line through $O$ in $\R^\ell$. The focal may have stable non isolated singularities (Cf. 
if $\Qbf:\R^2\to\R^4$). In Fig.\ref{fig:admissible-singul}, we show some singularities of hyperbolic surfaces (and a forbidden one). 

\begin{figure}[ht!]
	\centering
	\includegraphics [scale=0.2]{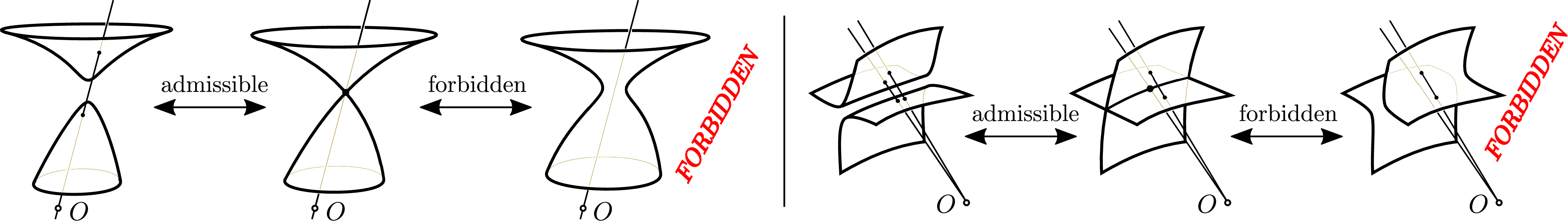}
	\caption{\small Admissible and forbidden local metamorphosis of hyperbolic singularities of the focal. }
	\label{fig:restricted-nbifurcation3D}
\end{figure}

If the singularities of $\mathcal{F}_{\Qbf}$ are unstable, 
the metamorphoses they may undergo are also restricted 
by the hyperbolicity of $\mathcal{F}_{\Qbf}$: 
given a line joining the origin to a hyperbolic singular point of $\mathcal{F}_{\Qbf}$ and 
taking a local smoothing  $\mathcal{F}_{\Qbf}(\varepsilon)$, the number of 
intersection points of every neighbouring line with $\mathcal{F}_{\Qbf}(\varepsilon)$, 
counting multiplicities, cannot change. So, those lines cut both local branches. This forbids some usual metamorphosis - see Fig.\,\ref{fig:restricted-nbifurcation3D}. 

\begin{figure}[ht!]
	\centering
	\includegraphics [scale=0.15]{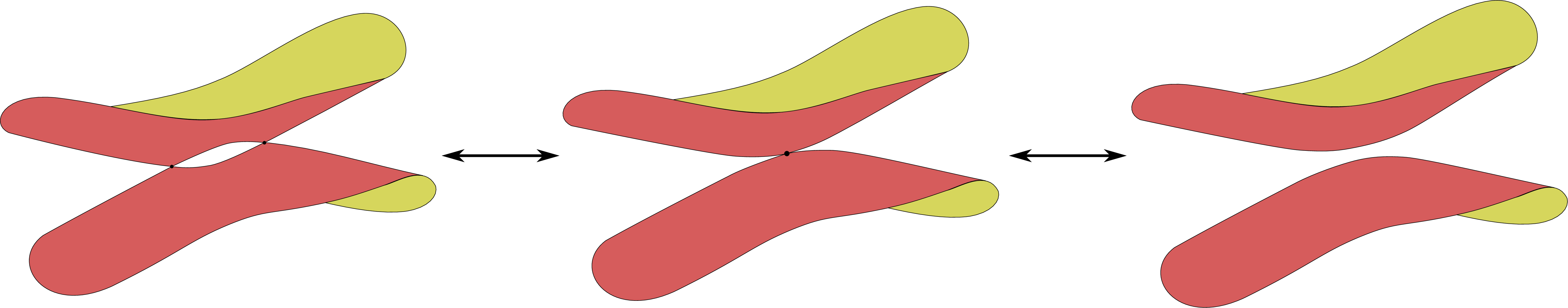}
	\caption{\small Creation/anihilation of two hyperbolic $A_1$ singularities of the focal at an $A_3$ singularity. }
	\label{fig:2A1_A3-bifurcation}
\end{figure}

\begin{example*}[Cayley Cubic Surface]
	Take the quadratic map $\Qbf:\R^3\to\R^3$, given by 
	\begin{equation}\label{eq:cubic-Steiner-Cayley}
		\Qbf(t)=(2t_2t_3, 2t_1t_3, 2t_1t_2)\,.
	\end{equation} 
	
	The focal $\mathcal{F}_{\Qbf}$ is given by the equation 
	$\begin{vmatrix}
		-1 &  z  & y \\ 
		z & -1  & x \\
		y &  x  & -1
	\end{vmatrix}=0\,,$
	that is, by the equation $x^2+y^2+z^2+2xyz-1=0$. It is the \textit{Cayley Cubic Surface} -  the unique cubic surface in $\RP^3$ 
	(up to projective transformation) which has four Morse conic singularities. 
	It is the union of one topological sphere with one topological projective plane which overlap at the four singular points. 
	\begin{figure}[ht!]
		\centering
		\includegraphics[width=0.42\textwidth]{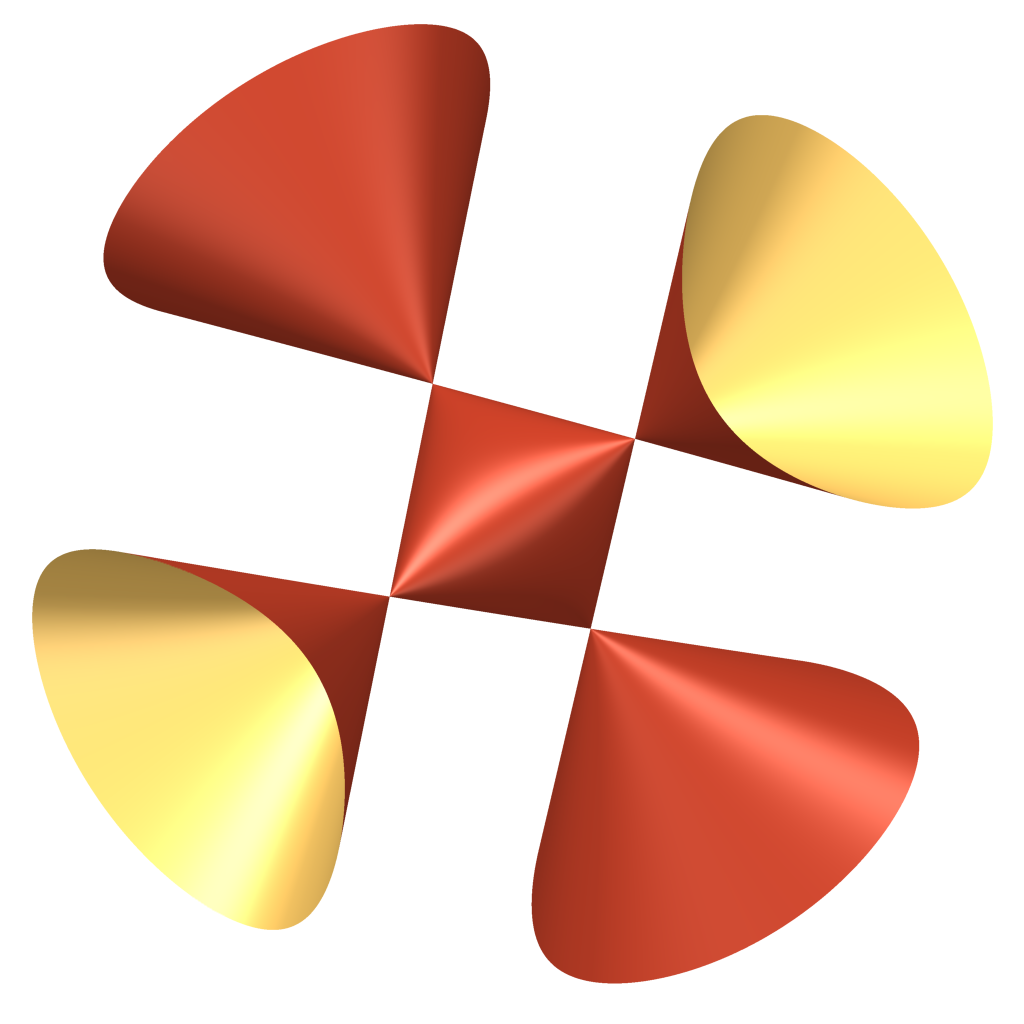}
		\caption{\small The Cayley Cubic Surface is the focal of the map $\Qbf(t)=(2t_2t_3, 2t_1t_3, 2t_1t_2)$.}
		\label{fig:Cayley-Cubic-Surface}
	\end{figure}
	
The four $A_1$-singular points of Cayley surface are stable under small deformations of $\Qbf$ 
in the space of quadratic maps - one can check it experimentally. For a quadratic map $\Qbf:\R^3\to\R^3$ we prove below that its focal is connected (implying that singular points are unavoidable) - moreover, 
the focals having two or four $A_1$ singular points are only stable ones. 

We have verified experimentally that under large generic deformations of $\Qbf$ two $A_1$-singular points 
may collapse at an $A_3$-singularity and disappear 
- Fig.\,\ref{fig:2A1_A3-bifurcation}. It is also possible 
that three $A_1$-singular points collapse at an $A_5$-singular point, but they don't disappear. 
\end{example*}

\subsection*{E \, Polar Dual Subvarieties} 
Steiner Roman Surface and Cayley Cubic Surface (the front $\Sigma_{\Qbf}$ and focal $\mathcal{F}_{\Qbf}$ of map 
\eqref{eq:cubic-Steiner-Cayley}) both with tetrahedral symmetry, have a geometric relation. 
Indeed, the front $\Sigma_{\Qbf}$ and focal $\mathcal{F}_{\Qbf}$ of every 
quadratic map $\Qbf:\R^k\to\R^\ell$ are related by projective duality. 

Each quadric $\mathcal{S}$ in $\RP^\ell$ provides a bijection between the hyperplanes in $\RP^\ell$ (the elements of $(\RP^\ell)^\vee$) and the points in $\RP^\ell$:   
given a point $P$ in $\RP^\ell$, a variable line $\lbf$ containing $P$ meets $\mathcal{S}$ at two points, say, $P_1(\lbf)$ and $P_2(\lbf)$. 
Write $\hat{P}(\lbf)$ for the point in $\lbf$ which is harmonic conjugate to $P$ with respect to $P_1(\lbf)$ and $P_2(\lbf)$. 
\textit{The geometric locus of the point $\hat{P}(\lbf)$ is a 
	hyperplane} (cf. \cite{Salmon}) - we note this hyperplane $H_P$. 
\smallskip

\noindent
\textit{\textbf{{\small Pole and Polar}}}.
The above hyperplane $H_P$ is called the \textit{polar} of the point $P$ with respect to $\mathcal{S}$. And $P$ is the \textit{pole} of the hyperplane $H_P$. 
\medskip

\noindent
\textit{\textbf{{\small Polar Duality}}}.
This bijection between the hyperplanes in $\RP^\ell$ and the points in the same 
space $\RP^\ell$, using a quadric $\mathcal{S}$, is an isomorphism called {\em polar duality}
\[\Psi_{\mathcal{S}} : (\RP^\ell)^\vee \stackrel{\simeq}{\longrightarrow}\RP^\ell\,.\]  

\begin{figure}[ht!]
	\centering
	\includegraphics[width=1\textwidth]{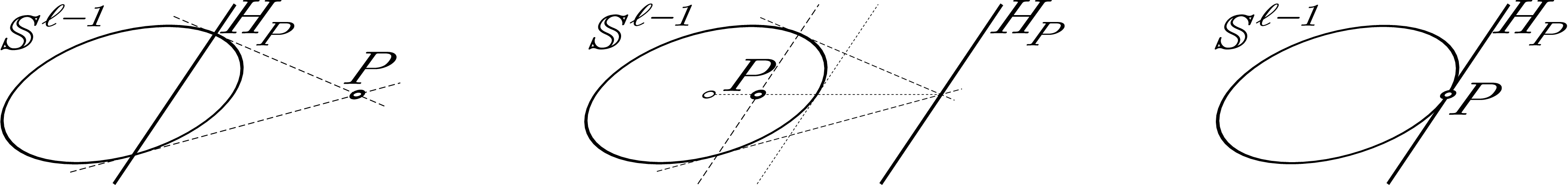}
	\caption{\small Correspondence Poles $\leftrightarrow$ Polars.}
	\label{fig:correspondence Poles-Polars}
\end{figure}

We take as $\mathcal{S}$ the unit hypersphere $\sph^{\ell-1}:=\{q\in\R^\ell:\langle q,q\rangle=1\}$ 
(Fig.\,\ref{fig:correspondence Poles-Polars}). 
In this case,      
\textit{the polar of the point $\,P=(p_1,\ldots,p_\ell)\neq O_{\R^\ell}\,$ has equation} 
$\,p_1q_1+\ldots+p_\ell q_\ell=1$:  
\begin{equation}\label{eq:pole-polar-relation}
	\mbox{(pole) }\, P \longleftrightarrow \mbox{(polar) }\, H_P -\mbox{hyperplane of equation} \,\, \langle P, q \rangle=1 \,. 
\end{equation} 

\begin{example*} 
	In $\R^3$, the polar plane of the point $P=(3,4,5)$ is given by the equation \,$3q_1+4q_2+5q_3=1$. 
\end{example*} 
%

\begin{remark}[Projective completion]\label{rem:projective-completion} 
Polar duality in relation \eqref{eq:pole-polar-relation} extends to $\RP^\ell\,(=\R^\ell\sqcup \RP^{\ell-1}_\infty)$: the poles of the hyperplanes through $O_{\R^\ell}$ (in $\R^\ell\sqcup \RP^{\ell-1}_\infty$) form
the hyperplane $\RP^{\ell-1}_\infty$ (at infinity), and the origin $O_{\R^\ell}$ is the pole of $\RP^{\ell-1}_\infty$. 
In homogeneous coordinates $[q_1:\ldots:q_\ell:q_{\ell+1}]$ the affine relation \eqref{eq:pole-polar-relation} becomes $\langle P,q\rangle=p_{\ell+1}q_{\ell+1}$ and the hyperplane $\RP_\infty^{\ell-1}$ has equation $q_{\ell+1}=0$.  
For any $q\in\R^\ell\sqcup\RP_\infty^{\ell-1}$, $q\neq O_{\R^\ell}$, the line joining $q$ to the origin $O_{\R^\ell}$
is orthogonal to the polar hyperplane of $q$. In particular, a hyperplane $H_q$ in $\R^\ell\sqcup\RP_\infty^{\ell-1}$ 
contains the origin of $\R^\ell$ iff its pole $q$ is at infinity in the direction orthogonal to $H_q$.  
\end{remark}

\noindent
\textit{\textbf{{\small Reciprocal Relation}}}. \textit{A point $R$ lies on the polar of $P$ iff $P$ lies on the polar of $R$}: 
\[R\in H_P\, \mbox{ iff }\, P\in H_R\,. \]
In consequence, \textit{the polar of $P$, $H_P$, consists of the poles of the hyperplanes through $P$}.
\medskip

\noindent
\textit{\textbf{{\small Polar Dual}}}.
The \textit{polar dual} of a subvariety $\Sigma$ in $\RP^\ell$ (with respect to $\mathcal{S}$) is the subvariety $\Sigma^\vee$ in $\RP^\ell$ 
which consists of the poles of all hyperplanes tangent to $\Sigma$: 
\[\Sigma^\vee=\left\{P\in\RP^\ell:H_P \mbox{ is tangent to } \Sigma \mbox{ at some } 
Q\in\Sigma \right\}\,. \] 

\begin{example*}[Quadrics]
It is classically known that the polar dual of a smooth quadric in $\RP^\ell$ is also a smooth quadric in $\RP^\ell$. 
\end{example*}

\noindent
\textbf{\small Fundamental Fact}: $\left(\Sigma^\vee\right)^\vee=\Sigma$\,. 

\section{Duality Theorem and Geometry of the Front and Focal}
\subsubsection*{A \ Duality Theorem}

\begin{theorem*}[Duality Theorem]
For any homogeneous quadratic map $\Qbf:\R^k\to\R^\ell$ the focal hypersurface $\mathcal{F}_{\Qbf}$ 
is the polar dual of the front $\Sigma_{\Qbf}:$ \,  
$\Sigma_{\Qbf}^\vee=\mathcal{F}_{\Qbf}$. 
\end{theorem*} 

\begin{corollary*}
The Cayley Cubic Surface is the polar dual of the Steiner Roman Surface. 
\end{corollary*}

\begin{proof}[\underline{Proof of Duality Theorem}]
	The polar hyperplane $H_{\Qp([t])}$ of the point $\Qp([t])$ is given by the equation 
	\begin{equation}\label{eq:polar-plane}
		\langle\,\Qp([t])\,,\,q\,\rangle=1\,. 
	\end{equation}
	But since $\Qbf$ is given in form \eqref{eq:homogeneous-Rk_Rn}, equation \eqref{eq:polar-plane} of the hyperplane $H_{\Qp([t])}$ is equivalent to 
	\[\sum_{i=1}^kt_i^2\langle\,\Abf_{ii}\,,\,q\,\rangle+
	\sum_{1\leq i<j\leq k}2\,t_it_j\langle\,\Abf_{ij}\,,\,q\,\rangle=\sum_{i=1}^kt_i^2, \ \mbox{ that is,}\]
	\begin{equation}\label{eq:polar-plane-f}
		\Phi_q(t):=\sum_{i=1}^kt_i^2\left(\langle\,\Abf_{ii}\,,\,q\,\rangle-1\right)+
		\sum_{1\leq i<j\leq k}2\,t_it_j\langle\,\Abf_{ij}\,,\,q\,\rangle=0
	\end{equation}
	
	This hyperplane $H_{\Qp([t])}$ is the union of the poles of all hyperplanes through $\Qp([t])$.  
	Namely, the coordinates $(q_1,\ldots,q_n)$ of a point $q$ satisfying \eqref{eq:polar-plane} 
	are the slope coefficients of its polar hyperplane $H_q$, which is a hyperplane 
	passing through $\Qp([t])$.
	
	To each point $\Qp([t])$ in $\Sigma_{\Qbf}$, polar duality associates the subset of the polar dual 
	$\Sigma_{\Qbf}^\vee$ formed by the poles of the hyperplanes that contain the tangent space to 
	$\Sigma_{\Qbf}$ at $\Qp([t])$. So this subset is an affine subspace $\Pi_{[t]}\Qbf$ of the polar hyperplane $H_{\Qp([t])}$.
	To get the equations of this subspace $\Pi_{[t]}\Qbf$, we note that 
	the tangent space $T_{\Qp([t])}\Sigma_{\Qbf}$ is generated by the partial derivative vectors $\D\Qbf/\D t_i$;  
	thus, besides eq.\,\eqref{eq:polar-plane-f}, that is $\Phi_q(t)=0$, we add $k-1$ equations from the condition that $k-1$ partial 
	derivatives of $\Phi_q(t)$ (with respect to $k-1$ of the $k$ variables $t_i$) vanish. 
	But these $k$ equations are equivalent to require that the $k$ partial derivatives of $\Phi_q(t)$ vanish, 
	due to the homogeneity of $\Phi_q(t)$ and Euler equation $\sum t_i(\D\Phi_q(t)/\D t_i)=2\Phi_q(t)$. 
	
	(Note that the subspace $\Pi_{[t]}\Qbf$ of $H_{\Qp([t])}$ is non empty iff $d_{[t]}\Qp$ is non surjective
	because $\Pi_{[t]}\Qbf$ is the dual of the image of the differential $d_{[t]}\Qp:T_{[t]}\RP^{k-1} \to T_{\Qp([t])}\R^\ell$.) 
	
	So taking the $k$ partial derivatives of \eqref{eq:polar-plane-f} and  dividing by $2$ we get: 
	\begin{alignat*}{5}
		t_1(\langle\Abf_{11},q\rangle-1)&+t_2\langle\Abf_{12},q\rangle&+t_3\langle\Abf_{13},q\rangle+&\ldots+t_k\langle\Abf_{1k},q\rangle &= 0\, \, \,\\
		t_1\langle\Abf_{21},q\rangle&+t_2(\langle\Abf_{22},q\rangle-1)&+t_3\langle\Abf_{23},q\rangle+&\ldots+t_k\langle\Abf_{2k},q\rangle &= 0\, \, \,\\
		\vdots & &\vdots & & \vdots \, \, \,\\
		t_1\langle\Abf_{k1},q\rangle&+t_2\langle\Abf_{k2},q\rangle&+t_3\langle\Abf_{k3},q\rangle+&\ldots+t_k(\langle\Abf_{kk},q\rangle-1) &\, \, \, =  0\,.
	\end{alignat*}
	This system of equations is equivalent to 
	\[
	\begin{pmatrix}
		\langle \Abf_{11}\,, \, q\rangle -1 & \langle \Abf_{12}\,, \, q\rangle & \ldots & \langle \Abf_{1k}\,, \, q\rangle\\ 
		\langle \Abf_{21}\,, \, q\rangle & \langle \Abf_{22}\,, \, q\rangle -1 & \ldots & \langle \Abf_{2k}\,, \, q\rangle\\
		\vdots & \vdots & \ddots & \vdots \\
		\langle \Abf_{k1}\,, \, q\rangle & \langle \Abf_{k2}\,, \, q\rangle & \ldots & \langle \Abf_{kk}\,, \, q\rangle-1
	\end{pmatrix}
	\begin{pmatrix}
		t_1 \\
		t_2 \\
		\vdots \\
		t_k
	\end{pmatrix}=
	\begin{pmatrix}
		0 \\
		0 \\
		\vdots \\
		0
	\end{pmatrix}\,,
	\]
which is an equation of the form $(M_{\Qbf}^q-\mathrm{Id}_k)t=\bar{0}$, where 
$(M_{\Qbf}^q)_{ij}=\langle\,\Abf_{ij}\,,\,q\,\rangle$. At the points $q\in\R^\ell$ where this 
equation holds (with $t\neq\bar{0}$), the matrix $\,M_{\Qbf}^q-\mathrm{Id}_k\,$ has non trivial kernel. 
This means that the polar dual of the front $\Sigma_{\Qbf}$ is the degree $k$ algebraic hypersurface 
of $\RP^\ell$ given (in affine coordinates) by the equation 
\[\det(M_{\Qbf}^q-\mathrm{Id}_k)=0,\]
which determines the focal $\mathcal{F}_{\Qbf}$ - see  \eqref{eq:focal_set-equation}.  
Duality Theorem is proved. 
\end{proof}

The duality works in both directions: $\Sigma_{\Qbf}^\vee=\mathcal{F}_{\Qbf}$ and $\mathcal{F}_{\Qbf}^\vee=\Sigma_{\Qbf}$. 
It holds also if the front $\Sigma_{\Qbf}$ or the focal $\mathcal{F}_{\Qbf}$ degenerate. 
Notice that they have different natures: the focal $\mathcal{F}_{\Qbf}$ is given by eq.\,\eqref{eq:focal_set-equation},  but the front $\Sigma_{\Qbf}$ is the image of a map $\RP^{k-1}\to\R^\ell$. 
\medskip

\noindent
\textbf{\small\textit{Tangent Hyperplanes to the Image of a Map}}. The \textit{tangent space} to the image of a smooth map $f:M^k\to\R^\ell$, $f(M)$, at $f(p)$ is the image of the differential of $f$ at $p$, 
$d_pf(T_pM)$, that is, $T_{f(p)}f(M):=d_pf(T_pM)$. \ 
An affine hyperplane $H$ in $\R^\ell$ is said to be \textit{tangent} to 
$f(M)\subset\R^n$ at $f(p)$ if it contains the tangent space to $f(M)$ at $f(p)$: $T_{f(p)}f(M)\subset H$. 

\begin{Example}\label{ex:ellipse-cone}
	The dual to an ellipse $\Qp:\RP^1\to\RP^2$ which degenerates into a double segment, consists of the two 
	lines whose points are the poles of all lines tangent to 
	that segment at its end-points (i.e. the polar lines of the end-points). 
	
	Conversely, observe that a focal $\mathcal{F}_{\Qbf}$ consisting of two lines is unstable but, by its hyperbolicity (Th.\,\ref{th:hyperbolicity-theorem}a), 
	it cannot undergo a Morse bifurcation (see Fig.\,\ref{fig:half-morse}). 
	\begin{figure}[ht!] 
		\centering
		\includegraphics[scale=0.17]{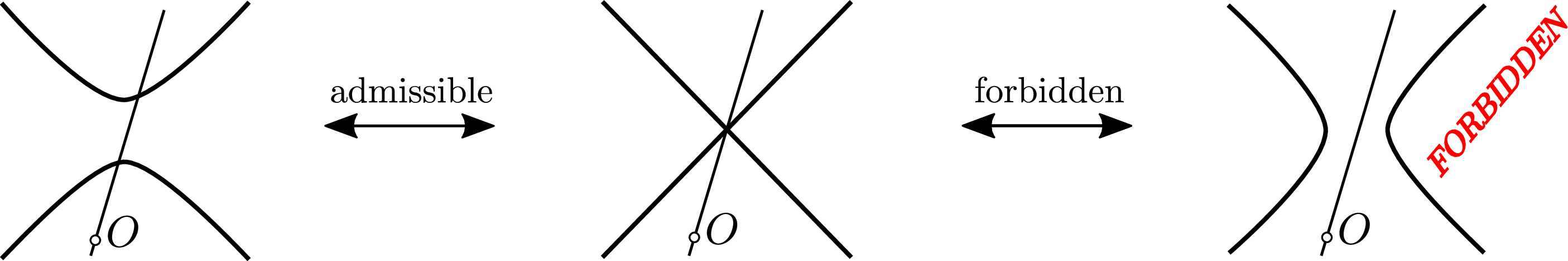}
		\caption{\small Permited and forbidden transitions of a degenerate conic $\mathcal{F}_{\Qbf}$ (by its hyperbolicity).}
		\label{fig:half-morse}
	\end{figure}

	To define its dual, we include $\mathcal{F}_{\Qbf}$ in a generic 
	$1$-parameter family of conics $\mathcal{F}_{\Qbf(\e)}$ (having the hyperbolicity property) and take the limit of their duals.  
	\begin{figure}[ht!] 
		\centering
		\includegraphics[scale=0.17]{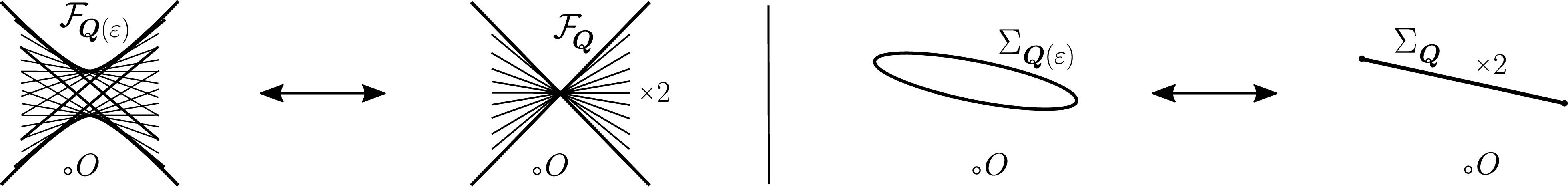}
		\caption{\small Permited and forbidden transitions of a degenerate conic $\mathcal{F}_{\Qbf}$ (by its hyperbolicity).}
		\label{fig:pensil-bufurcation}
	\end{figure}
	At the limit, the tangent lines to the hyperbola $\mathcal{F}_{\Qbf(\e)}$ become a segment of the pencil of lines through the double point, 
	taken with multiplicity two (one for each branch) - Fig.\,\ref{fig:pensil-bufurcation}. 
	The poles of those lines form a segment in the polar line of the double point, taken with multiplicity two. 
\end{Example}

\begin{Example}[Cone-Ellipsoid]\label{ex:ellipse-cone}
The dual to an ellipse $\Sigma$ in $\RP^3$ is a non-degenerate quadratic cone $\mathcal{F}$ whose vertex 
is the pole of the plane containing $\Sigma$. 
Conversely, the poles of the planes tangent to the cone $\mathcal{F}$ belong to the polar hyperplane of its 
vertex, where they form the original ellipse $\Sigma$. 
Similarly, the dual to a quadratic cone $\mathcal{F}$ of type $(+,\ldots,+,-)$ in $\R^\ell$ is a codimension $2$ 
ellipsoid lying in the polar of the vertex of $\mathcal{F}$. 
	
For example, the dual to the cone of equation 
$a_1^2q_1^2+\cdots+a_{\ell-1}^2q_{\ell-1}^2-(q_\ell-r)^2=0$
is the codimension $2$ ellipsoid given by 
the following two equations: 
\[\frac{q_1^2}{a_1^2}+\cdots+\frac{q_{\ell-1}^2}{a_{\ell-1}^2}=\frac{1}{r^2}, \qquad q_\ell=\frac{1}{r}\,.\]
\end{Example}

\noindent
\textbf{\small Euclidean-Projective Remark}. 
The duality relation between $\mathcal{F}_{\Qbf}$ and $\Sigma_{\Qbf}$ is not only projective. 
The Euclidean structure of $\R^\ell$ enters in the definitions of 
$\Sigma_{\Qbf}$, of $\mathcal{F}_{\Qbf}$ 
and of the polar duality with respect to the unit hypersphere. 
This Euclidean dependence makes rigid the duality relation between 
$\mathcal{F}_{\Qbf}$ and $\Sigma_{\Qbf}$: each point of one of them is the pole of a specific hyperplane in Euclidean space $\R^\ell$ which is tangent to the other. 
\medskip

Our proof of Duality Theorem ``excludes'' the case in which $\Qp([t])=O_{\R^\ell}$ because we give it using affine coordinates. 
To include this case, it suffices to use homogeneous coordinates (see Remark\,\ref{rem:projective-completion}). 

\begin{corollary*}[of Duality Theorem]\label{th:dual-Veronese}
	The polar dual of a Veronese $(k-1)$-manifold is a degree $k$ hyperbolic hypersurface. 
	It is the union of a $(k-1)$-parameter family of 
	$\frac{1}{2}k(k-1)$-dimensional projective subspaces that have one common point. 
\end{corollary*} 

\begin{Example}\label{ex:Veronese-Manifolds-2-3}
	$\boldsymbol{(k=2)}$ A Veronese $1$-manifold is an ellipse in $\R^3$ (the image of $\Qp:\RP^1\to\R^3$). Its dual $\mathcal{F}_{\Qbf}$ is quadratic cone in $\RP^3$ (a $1$-parameter family of lines). \\
	$\boldsymbol{(k=3)}$ A Veronese $2$-manifold is the image of a quadratic map $\Qp:\RP^2\to\R^6$. 
	Its dual $\mathcal{F}_{\Qbf}$ is a cubic hypersurface which
	is a union of a $2$-parameter family of $3$-dimensional projective subspaces 
	which have one common point. 
\end{Example}  

\subsubsection*{B \ Unavoidable Singularities of the Focal $\mathcal{F}_{\Qbf}$ if $k\leq \ell$ }\label{sect:tot-reducible}
If the focal $\mathcal{F}_{\Qbf}$ is smooth, then 
its nested hyperpheres (and the hyperplane
for odd $k$) are disjoint. If $\mathcal{F}_{\Qbf}$ is non smooth, its singular points 
are common points of consecutive nested hyperspheres. 
If $k\leq \ell$, some additional restrictions are imposed to $\mathcal{F}_{\Qbf}$. 

\begin{theorem}\label{th:Unavoidable-singularities-FQ}
Let $\Qbf:\R^k\to\R^\ell$ a homogeneous quadratic map with $k\leq \ell$. 
	
\noindent
$(a)$ For each $i\in\N$ with $i\leq\lfloor\frac{k-1}{2}\rfloor$, the nested topological hyperspheres $\phi_i$ and $\phi_{i+1}$ in the focal $\mathcal{F}_{\Qbf}$ 
(and also the topological hyperplane $\eta$ if $k$ is odd) overlap at singular points.
\smallskip

\noindent
$(b)$ If $k<\ell$ and $\Qp(\RP^{k-1})$ is not contained in a $(k-1)$-affine subspace of $\R^\ell$, then $\mathcal{F}_{\Qbf}$ is the union of a $(k-1)$-parameter family of $(\ell-k)$-dimensional projective subspaces in $\RP^\ell(=\R^\ell\sqcup\RP^{\ell-1}_\infty)$. 
\smallskip
	
\noindent
$(c)$ If $\Qbf:\R^k\to\R^k$ is in general position, then the focal $\mathcal{F}_{\Qbf}$  contains $\lfloor\frac{k-1}{2}\rfloor$ collections of 
locally hyperbolic singular points, such that $\phi_i$ and $\phi_{i+1}$ overlap 
at the points of the $i$-th collection $\{P^i_1,\ldots,P^i_{n_i}\}_{n_i\geq 1}$.
If $k$ is odd, the topological hypersphere 
$\phi_{\lfloor\frac{k-1}{2}\rfloor}$and the topological hyperplane $\eta$ 
also overlap at singular points.  
\end{theorem} 
\begin{proof}
$(a)$ If $k\leq \ell$, $\Sigma_{\Qbf}=\Qp(\RP^{k-1})$ is  connected. Hence, its dual $\mathcal{F}_{\Qbf}$ is connected, and then, any two consecutive nested hyperspheres 
of $\mathcal{F}_{\Qbf}$ (the last one with $\eta$ for odd $k$) have common singular points. 
\smallskip 

\noindent
$(b)$ If $\Qp(\RP^{k-1})$ is not contained in a $(k-1)$-dimensional affine subspace, then $\Sigma_{\Qbf}$ coincides with $\Qbf(\RP^{k-1})$, 
which is a $(k-1)$-dimensional subvariety in $\R^\ell$. At each point $p$ in $\Sigma_{\Qbf}$,
the bundle of hyperplanes that are tangent to $\Sigma_{\Qbf}$ at $p$ (i.e. that contain $T_p\Sigma_{\Qbf}$) is of dimension $n-(k-1)-1$. 
Thus the poles of those tangent hyperplanes form a $(n-k)$-dimensional projective subspace in $\R^\ell\sqcup\RP^{\ell-1}_\infty$. 
Therefore, $\mathcal{F}_{\Qbf}$ is the union of a $(k-1)$-parameter family of $(n-k)$-dimensional projective subspaces. 
\smallskip

\noindent
$(c)$ The focal is connected (because $\Sigma_{\Qbf}=\Qp(\RP^{k-1})$ is connected), 
but is not a union of projective subspaces and has no other restriction besides its hyperbolicity.  
Thus, by Theorem \ref{th:hyperbolic-singularities}, two consecutive nested topological hyperspheres overlap at 
a non-empty collection of stable isolated hyperbolic singular points {\rm (the Cayley Cubic surface, 
which is the focal of the map, has four such singular points and the focal $\mathcal{F}_{\Qbf}(\varepsilon)$ 
of a small deformation of $\Qbf(\varepsilon)$ is projectively equivalent to the Cayley surface)}. 
\end{proof}	
	
\subsubsection*{C \ Configuration of the Front $\Sigma_{\Qbf}$ and Toeplitz-Haussdorf Theorem}
\begin{figure}[ht!]
\centering
\includegraphics [scale=0.18]{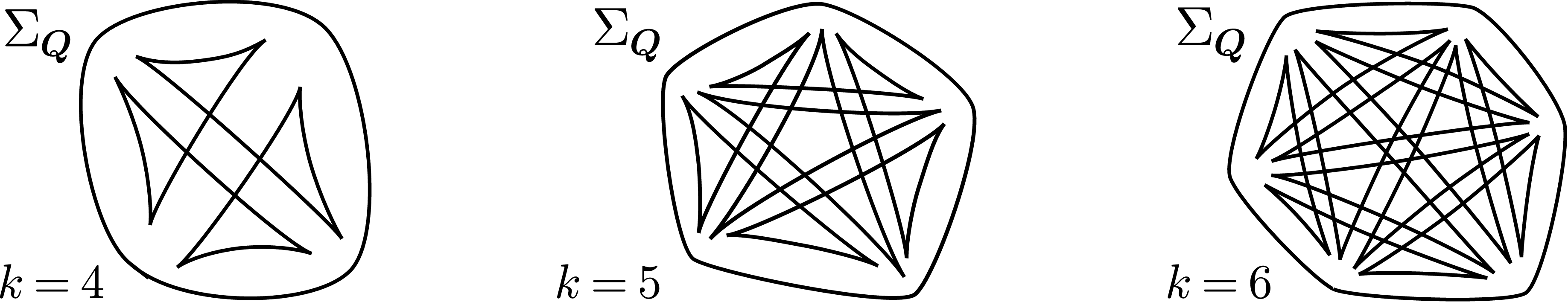}
\caption{\small Fronts 
$\Sigma_{\Qbf}$ of quadratic maps $\Qbf:\R^k\to\R^2$ for $k=4, 5, 6$, duals to the focals of Fig.\ref{fig:nested-smooth-spheres}.}
\label{fig:fronts_k=4-5-6}
\end{figure}

The dual version of Theorem\,\ref{th:topological-configurations} provides the configuration of the front $\Sigma_{\Qbf}$.  
\begin{theorem}\label{th:configuration-SigmaQ}
Let $\Qbf:\R^k\to\R^\ell$ a homogeneous quadratic map with $k>\ell$. 
\smallskip

\noindent
\textit{\textbf{a}}$)$ The front $\Sigma_{\Qbf}$ consists of $\lfloor\frac{k+1}{2}\rfloor$ immersed hypersurfaces 
{\rm  (figs.\,\ref{fig:fronts_k=4-5-6} and \ref{fig:fronts_k=5-5-6})}. \ 
\smallskip
	
\noindent
\textit{\textbf{b}}$)$ {\rm (Toeplitz-Haussdorf Theorem)}
The range of $\Qp:\RP^{k-1}\to\R^\ell$ is a convex subset in $\R^\ell$ whose boundary is the topological hypersphere $\sigma_1$ (the polar dual of $\phi_1$). 
\end{theorem}

\begin{proof}
\noindent
\textit{\textbf{a}}) 
The hypersurfaces forming $\Sigma_{\Qbf}$ are just the polar duals of the respective  $k/2$ (or $(k+1)/2$)
hypersurfaces $\phi_1$, $\phi_2$,\ldots, of $\mathcal{F}_{\Qbf}$. We note them by $\sigma_1$, $\sigma_2$, etc. 
\smallskip 

\noindent
\textit{\textbf{b}}) 
The boundary of the range of $\Qp$ is the outermost component $\sigma_1$ of the critical value set $\Sigma_{\Qbf}$ - it is the polar dual 
of the innermost component $\phi_1$ of $\mathcal{F}_{\Qbf}$. 
The topological hypersphere $\sigma_1$ is convex because $\phi_1$ is convex and duality respects convexity: 
\textit{if a hypersurface in $\RP^\ell$ bounds a convex set, then its dual hypersurface also bounds a convex set}.  
\end{proof}

\begin{proposition}
	If $k\leq \ell$ and $\Qbf:\R^k\to\R^\ell$ is a homogeneous quadratic map in general position, then the range $\Qp(\RP^{k-1})$ is not convex. 
\end{proposition}

\begin{proof}
	A necessary condition to be convex for a $(k-1)$-dimensional variety is to lie 
	in a $(k-1)$-affine subspace. But if $k-1<\ell$ and $\Qbf:\R^k\to\R^\ell$ is in 
	general position, the range $\Qp(\RP^{k-1})$ is not contained in a $(k-1)$-affine subspace. 	
\end{proof} 

\begin{Example}
For $\Qbf:\R^2\to\R^{\ell\geq 2}$ in general position, the range $\Sigma_{\Qbf}=\Qp(\RP^1)$ is a non degenerate ellipse: 
it is not contained in a line and is not a convex set. 
\end{Example}

\begin{corollary*}[of Theorem\,\ref{th:configuration-SigmaQ}]\label{Cor:range=closed-ball}
Let $\Qbf:\R^k\to\R^\ell$ a quadratic map with $k> \ell$. 
If the innermost topological hypersphere $\phi_1$ in the focal $\mathcal{F}_{\Qbf}$ is smooth, then 
the range of $\Qp$ $(\Qp(\RP^{k-1}))$ is a closed ball whose boundary is a strictly convex smooth ovaloid. 
\end{corollary*}

\begin{figure}[ht!]
\centering
\includegraphics [scale=0.2]{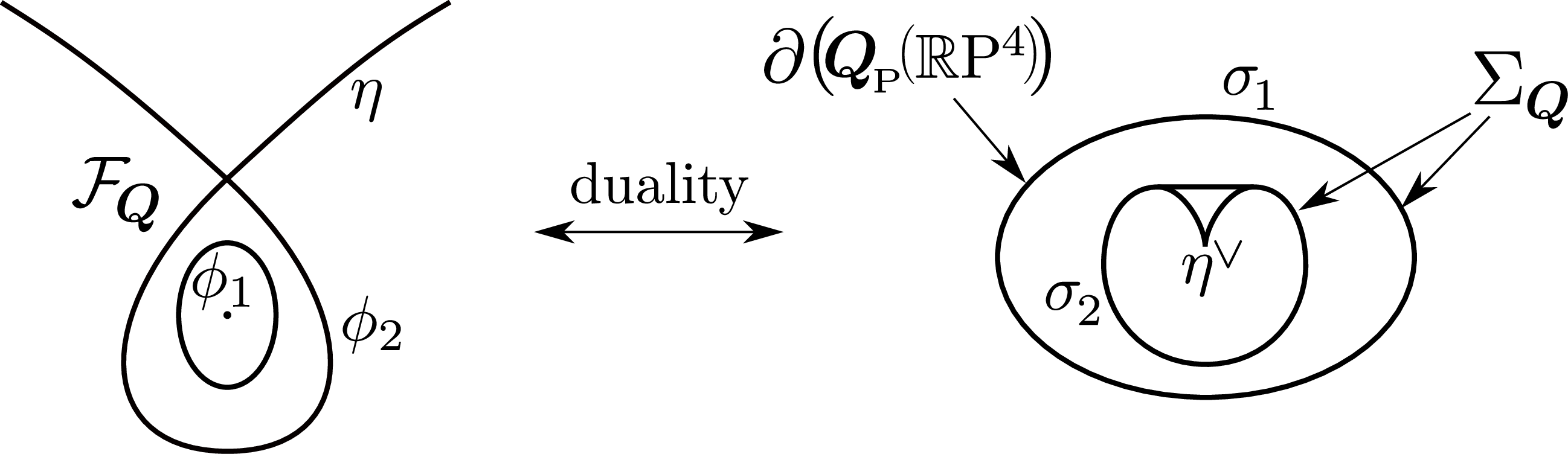}
\caption{\small The innermost component $\phi_1$ of $\mathcal{F}_{\Qbf}$ and its dual $\sigma_1$ in $\Sigma_{\Qbf}$, both are ellipses.}
	\label{fig:special-quintic}
\end{figure}

\begin{Example}\label{ex:special-quintic}
The boundary of the range of $\Qp$ for the quadratic map $\Qbf:\R^5\to\R^2$, 
\[\Qbf(t)=\left(\,\frac{t_1^2}{4}+\frac{t_2^2}{4}-\frac{t_3^2}{4}-\frac{t_4^2}{3}+\frac{t_5^2}{3}+\frac{2}{3}t_4t_5\, \, ; \,
-\frac{t_4^2}{2}+\frac{t_5^2}{2}+\frac{t_1t_2}{2}+\frac{t_2t_3}{2}-t_4t_5\,\right)
\]
is an ellipse - the outermost component of the front $\Sigma_{\Qbf}$ (Fig.\,\ref{fig:special-quintic}-right).  
Its focal $\mathcal{F}_{\Qbf}$ is the singular quintic of Fig.\,\ref{fig:special-quintic}-left. 
\end{Example}

\subsubsection*{D \ Polytopes (Fronts $\Sigma_{\Qbf}$) and Hyperplane Arrangements (Focals $\mathcal{F}_{\Qbf}$)}\label{sect:tot-reducible}

For quadratic maps $\Qbf:\R^k\to\R^\ell$, the so called \textit{diagonal case} is one of the most 
studied cases because it is algebraically more simple (but being very degenerate).

\begin{theorem}\label{th:dual-totally-reducible}
Let $\Qbf:\R^k\to\R^\ell$ a quadratic map with $k> \ell$. 
If the focal $\mathcal{F}_{\Qbf}$ is the union of $k$ hyperplanes in general position, 
then the front $\Sigma_{\Qbf}$ consists of all the facets {\rm (of dimension $d<\ell$)} of the polytope whose vertices are the poles of those $k$ hyperplanes. 
Each $d$-facet appears with multiplicity $2^d$ $($or with multiplicity $2^{d+1}$ for $\Qs:\sph^{k-1}\to\R^\ell$$)$. 
\end{theorem} 

\begin{proof}
Indeed, the poles of the $k$ hyperplanes that form $\mathcal{F}_{\Qbf}$ are $k$ points in $\R^\ell$ in general position 
(no $\ell+1$ points are contained in a hyperplane). Since no hyperplane of $\mathcal{F}_{\Qbf}$ contains the origin, 
we can always write the equation of the $i$-th hyperplane of $\mathcal{F}_{\Qbf}$ as $a_i^1x_1+a_i^2x_2+\ldots+a_i^\ell x_\ell=1$;
thus, its pole is $\Abf_i=(a_i^1,a_i^2,\ldots,a_i^\ell)$ and the quadratic map is written as
\begin{equation}\label{eq:diagonal-quad-map}
\Qbf(t)=\sum_{i=1}^kt_i^2\Abf_{i}, 
\quad \mbox{where } \|t\|^2=\sum_{i=1}^kt_i^2\,,
\end{equation} 
which is the so called \textit{diagonal case}. 
	
The differential of $\Qs:\sph^{k-1}\to\R^\ell$ at $t\in\sph^{k-1}$, $d_t\Qs: T_t\sph^{k-1}\to T_{\Qs(t)}\R^\ell$, is given by 
\begin{equation}\label{eq:differential-diag}
d_t\Qs(v)=2\Sigma_{i=1}^kt_i\Abf_idt_i(v)=2\Sigma_{i=1}^kt_i\Abf_idt_iv_i, 
\end{equation}
where $t\in\sph^{k-1}\subset\R^k$, $v\in T_t\sph^{k-1}\simeq\R^{k-1}$ (which imply $\langle\,t\,,\,v\,\rangle=0$). 
	
Since the points $\Abf_i$ are in general position, the differential\,\eqref{eq:differential-diag} at $t$ is not 
surjective iff the number of non zero components of $t=(t_1,\ldots,t_k)$ is less than $\ell$, that is, 
iff $t$ belongs to a coordinate subspace in $\R^k$ given by 
$t_{i_1}=0$, $t_{i_2}=0$, \ldots, $t_{i_m}=0$, where $m\geq k-\ell+1$. 
Such a subspace intersects $\sph^{k-1}$ along a sphere $\sph^d$ (with $d=k-m-1$) 
whose image $\Qs(\sph^d)$ is a $d$-facet $\{t_{j_1}^2\Abf_{j_1}+t_{j_2}^2\Abf_{j_2}+\ldots+t_{j_{d+1}}^2\Abf_{j_{d+1}}\}$  
of the polytope $\Qs(\sph^{k-1})$. 
	
A point $p\in\R^\ell$ belongs to such an $d$-facet iff there exist $t\in\R^k$, with exactly $d+1$ non zero components 
$t_{j_1},t_{j_2},\ldots,t_{j_{d+1}}$ with $t_{j_1}^2+t_{j_2}^2+\ldots+t_{j_{d+1}}^2=1$ 
and such that $t_{j_1}^2\Abf_{j_1}+t_{j_2}^2\Abf_{j_2}+\ldots+t_{j_{d+1}}^2\Abf_{j_{d+1}}=p$. 
In such case, the point $p$ is a critical value for $2^{d+1}$ distinct critical points $t$ in $\sph^{k-1}\subset\R^k$. 
Namely, for each choice of the signs $\pm t_{j_i}$ for the $d+1$ non-zero components of $t$ 
(among the $2^{d+1}$ possible choices), we get a point $\bar{t}$ satisfying $\Qs(\bar{t})=p$ because $t_{j_i}^2=(-t_{j_i})^2$. 
If we replace $\Qs$ with the map $\Qp:\RP^{k-1}\to\R^\ell$, then $p$ is a critical value for $2^d$ critical points in $\RP^{k-1}$. 
	
This means that each $d$-facet of the polytope $\Qp(\RP^{k-1})$ must be taken with multiplicity $2^d$ 
for $\Qp:\RP^{k-1}\to\R^\ell$ (or with multiplicity $2^{d+1}$ for $\Qs:\sph^{k-1}\to\R^\ell$).  
\end{proof}

The polytopes of Theorem\,\ref{th:dual-totally-reducible} have been widely
studied as the bifurcation set of the topological type of the fibres $\Qs^{-1}(w)$, with $w\in\R^\ell$ (cf. \cite{Santiago1}).

\subsubsection*{E \ Class: The Front $\Sigma_{\Qbf}$ of a Quadratic Map $\Qbf:\R^k\to\R^\ell$ is of Class $k$ }\label{sect:tot-reducible}

\noindent
\textbf{\textit{{\small Class}}}. Let $M$ be an algebraic subvariety in $\RP^\ell$ whose dual $M^\vee$ is a hypersurface, and $\lbf^{\ell-2}$ be a codimension-$2$ 
projective subspace not tangent to $M$. The \textit{class} of $M$ is the number of real and imaginary hyperplanes that 
contain $\lbf^{\ell-2}$ and are tangent to $M$. (\textit{The class of $M$ does not depend on $\lbf^{\ell-2}$} and is equal to the degree of $M^\vee$.) 
\smallskip 

\noindent
\textbf{\textit{{\small Exterior Subspace}}}. Given a subvariety $M$ in $\R^\ell$, a codimension $2$ affine 
subspace $\lbf^{\ell-2}$, is said to be \textit{exterior} to $M$ 
iff there exist a hyperplane that contains $\lbf^{\ell-2}$ and that does not intersect $M$. 

\begin{example}
Given a nondegenerate ellipse $\Sigma$ in $\R^3$, a line $\lbf$ is exterior to $\Sigma$ iff $\lbf$ does not intersect the plane disc bounded by $\Sigma$.  
\end{example}

\begin{figure}[ht!]
\centering
\includegraphics [scale=0.075]{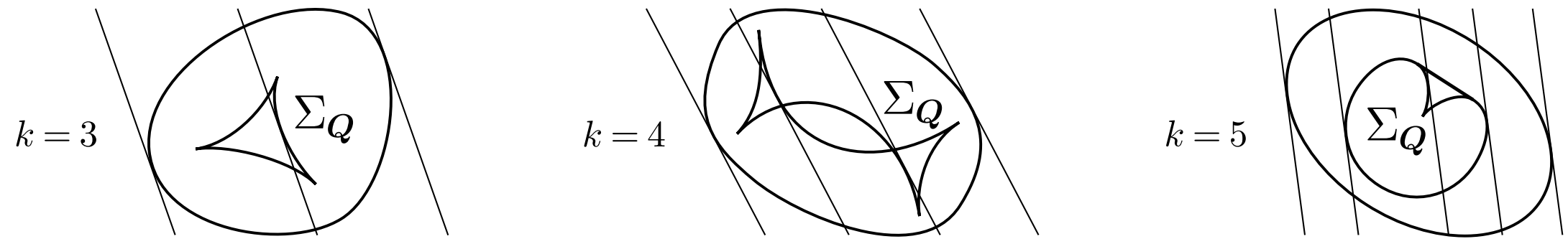}
\caption{\small The front of a quadratic map $\Qbf:\R^k\to\R^2$ is of class $k$. Here $k=3, 4, 5$.}
\label{fig:class-fronts}
\end{figure}

The dual version of Theorem\,\ref{th:hyperbolicity-theorem} (hyperbolicity of  $\mathcal{F}_{\Qbf}$) provides the class of $\Sigma_{\Qbf}$:

\begin{theorem}\label{th:class-k}
Let $\Qbf:\R^k\to\R^\ell$ be a homogeneous quadratic map. 
The front $\Sigma_{\Qbf}$ is an algebraic variety of class $k$.  
Moreover, through each $(\ell-2)$-subspace exterior to $\Sigma_{\Qbf}$ there pass $k$ real hyperplanes tangent 
to the front $\Sigma_{\Qbf}$ {\rm (Fig.\,\ref{fig:class-fronts})}. 
\end{theorem}

\begin{proof}
On one hand a point belongs to the open convex set bounded by the innermost topological hypersurface $\phi_1$ 
of $\mathcal{F}_{\Qbf}$ iff its polar hyperplane does not intersect $\Sigma_{\Qbf}$. 
On the other hand, the dual of a line $\lbf$ in $\RP^\ell$ is an $(\ell-2)$-dimensional subspace $L^{\ell-2}$ which consists of the 
poles of all hyperplanes containing $\lbf$, and, conversely, $\lbf$ consists of the poles of the hyperplanes containing $L$. 
Thus an $(\ell-2)$-dimensional subspace $L^{\ell-2}$ is exterior to $\Sigma_{\Qbf}$ iff 
its dual line $\lbf$ in $\RP^\ell$ contains a point of the set bounded by the innermost topological hypersurface $\phi_1$.  
Theorem\,\ref{th:topological-configurations}-\textit{\textbf{b}} states that 
the focal $\mathcal{F}_{\Qbf}$ is hyperbolic relative to every such point. 

Therefore $\lbf$ intersects $\mathcal{F}_{\Qbf}$ at $k$ real points that  are the poles of $k$ real hyperplanes through $L$, which are tangent to $\Sigma_{\Qbf}$. 
So Theorem\,\ref{th:class-k} is the dual formulation of Th.\,\ref{th:topological-configurations}-\textit{\textbf{b}}. 
\end{proof}


\begin{corollary}\label{cor:class-Veronese}
	Every Veronese submanifold $V_{k-1}$ is of class $k$, and through every codimension $2$ affine subspace 
	$($in $\R^{\frac{1}{2}k(k+1)}$$)$ exterior to $V_{k-1}$, in general position, there pass exactly $k$ real hyperplanes tangent to $V_{k-1}$.  	
\end{corollary}

\section{Duality, Parabolic Set and Front Singularities}
\subsubsection*{A \ Duality and Singularities of the Fronts}

\begin{figure}[ht!]
\centering
\includegraphics [scale=0.45]{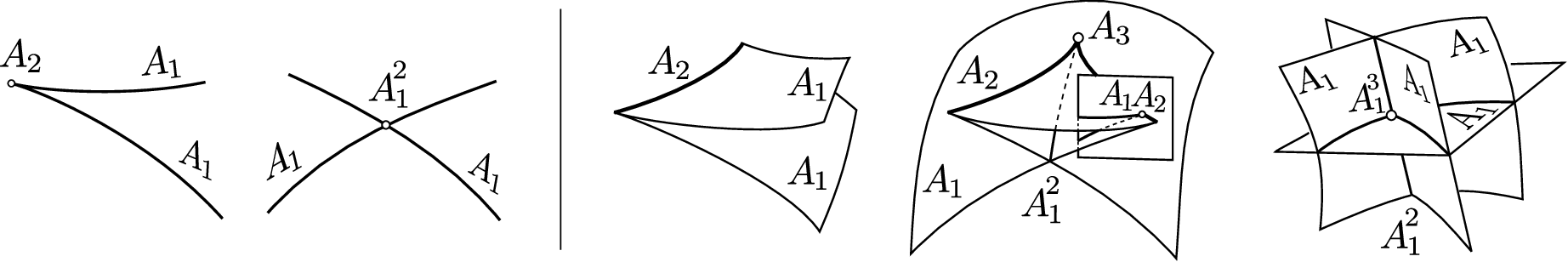}
\caption{\small Singularities of generic fronts in $2$-space and in $3$-space.}
\label{fig:front-singularities}
\end{figure}

The dual of an immersed submanifold in $\RP^n$ is the \textit{front} of a Legendrian map (cf. \cite{Arnoldwfeeml,avg}). 
The stable singularities of $1$-dimensional fronts in the plane are the semicubic cusps ($A_2$) and 
\textit{crossing points} ($A_1^2$) - 
Fig.\,\ref{fig:front-singularities}-left. 
The stable singularities of fronts in $3$-space are the semicubic cuspidal edge ($A_2$), the 
{\it swallowtail point} ($A_3$) and the transverse intersections of two or 
three branches ($A_1^2$), ($A_1A_2$), ($A_1^3$) - Fig.\,\ref{fig:front-singularities}-right. 

\textit{The stable singularities of fronts in $n$-space $(n\leq 6)$ are locally diffeomorphic to the discriminants of the Euclidean reflection groups $A_n$, $D_n$, $E_n$}
(cf. \cite{avg}). 

\subsubsection*{B \ Parabolic Set of the Focal $\mathcal{F}_{\Qbf}$ $\longleftrightarrow$ Singularities of the Front $\Sigma_{\Qbf}$}

\begin{figure}[ht!]
	\centering
	\includegraphics [scale=0.2]{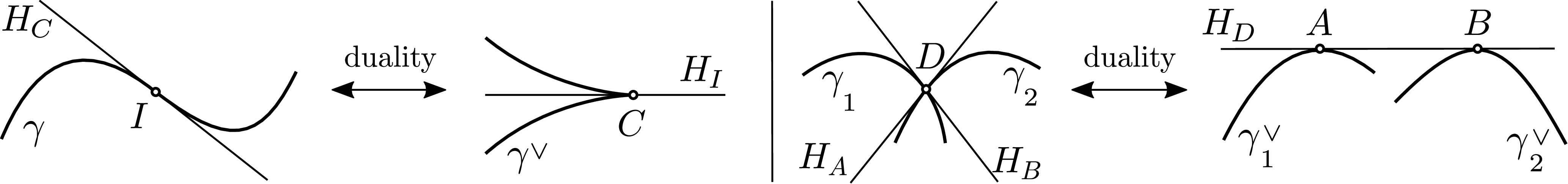}
	\caption{\small Duality correspondences: inflection $\leftrightarrow$ cusp \, and \, bitangent line $\leftrightarrow$ double point.}
	\label{fig:dual-correspondences}
\end{figure}

It is known (cf.\,\cite{Klein-elementary3}) that the tangent line $H_C$ to 
a plane curve $\g$ at an inflection $I$ is represented by a cusp $C$ in the dual curve $\g^\vee$, and vice-versa (Fig.\,\ref{fig:dual-correspondences}-left); and  a crossing $D$ of two branches of 
$\g$ represents a line $H_D$ 
which is tangent to the dual $\g^\vee$ at two points $A$, $B$,  and these points represent the lines $H_A$, $H_B$ tangent to $\g$ at $D$ 
(Fig.\,\ref{fig:dual-correspondences}-right).

\begin{Example}
	Take a homogeneous quadratic map $\Qbf:\R^3\to\R^2$. 
	By Theorem\,\ref{th:Focal-Eq}, the focal $\mathcal{F}_{\Qbf}$ is a curve of degree $3$. If $\mathcal{F}_{\Qbf}$ is smooth 
	and non reducible, by Th.\,\ref{th:topological-configurations}, its only possible 
	configuration in $\RP^2$ consists of an oval that contains the origin of $\R^2\subset\R^2\sqcup\RP^1_\infty$ 
	and a topological line with $3$ inflections (Fig.\,\ref{fig:dual-cubic}-left). 
	\begin{figure}[ht!]
		\centering
		\includegraphics [scale=0.075]{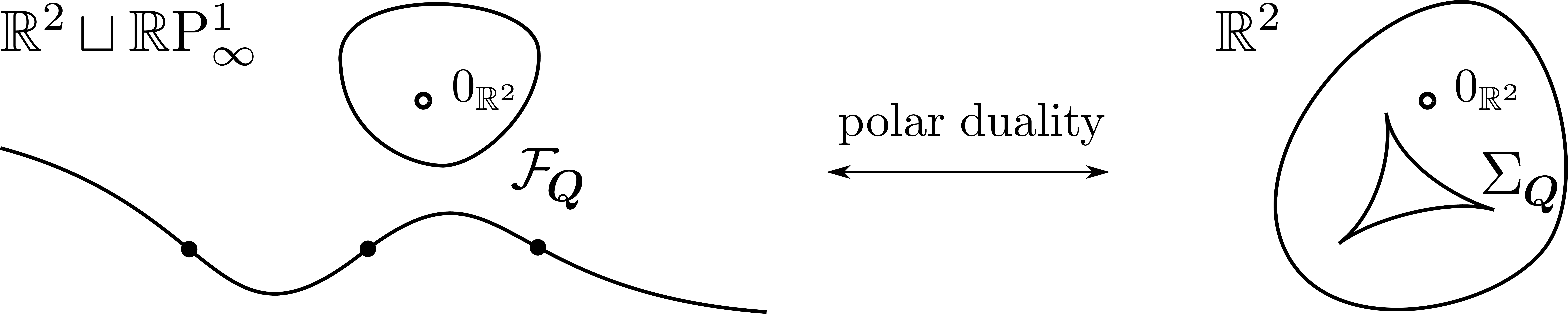}
		\caption{\small A cubic focal curve $\mathcal{F}_{\Qbf}\subset\R^2\sqcup\RP^1_\infty$ and 
			its polar dual $\mathcal{F}_{\Qbf}^\vee=\Sigma_{\Qbf}\subset\R^2$.}
		\label{fig:dual-cubic}
	\end{figure}
	Thus the front $\Sigma_{\Qbf}$ consists of an oval and a deltoid curve inside the oval, 
	whose $3$ cusps correspond to the $3$ inflections of the topological line of $\mathcal{F}_{\Qbf}$ (Fig.\,\ref{fig:dual-cubic}-right). 
\end{Example}

Similarly, \textit{the parabolic curve of a surface $S$ corresponds to the cuspidal edge of its dual surface $S^\vee$} (cf. \cite{Salmon}).
The notions of inflection point of a plane curve and of parabolic curve of a surface generalise to hypersurfaces in higher dimensional spaces: 
\medskip

\noindent
\textbf{\textit{{\small Parabolic Set}}}. 
The \textit{parabolic set} of a hypersurface $M$ in $\RP^\ell$ (or in $\R^\ell$) consists of the points $p$ in $M$ at which the tangent hyperplane has higher contact than usual with $M$. 
(This happens iff the second fundamental form of $M$ at $p$ is degenerate.) 

\section{Isolated singularities in $\mathcal{F}_{\Qbf}$ and $(\ell-2)$-ellipsoids in $\Sigma_{\Qbf}$} 

\subsection*{A \ $A_1$-singularities in $\mathcal{F}_{\Qbf}$ correspond to $(\ell-2)$-ellipsoids in $\Sigma_{\Qbf}$ ($k=\ell$)} 

The subvarieties $\Sigma_{\Qbf}$ and $\mathcal{F}_{\Qbf}$ in $\RP^\ell$ are dual to each other, but they have not similar properties. 
As a rule, the front $\Sigma_{\Qbf}$ has generic Legendre front singularities, while the focal 
$\mathcal{F}_{\Qbf}$ does not admit such singularities, due to its hyperbolicity 
(Th.\,\ref{th:hyperbolic-singularities}).  

Using example \ref{ex:ellipse-cone}, the dual version of Theorem\,\ref{th:Unavoidable-singularities-FQ}-$c$ provides 

\begin{theorem}\label{th:k-k-ellipsoids}
Let $\Qbf:\R^k\to\R^k$ be a homogeneous quadratic map in general position. 
\smallskip

\noindent
The front $\Sigma_{\Qbf}$ contains $\lfloor\frac{k-1}{2}\rfloor$ collections of 
ellipsoids of codimension $2$, which consist of infinitely degenerate parabolic points.
The ellipsoids of the $i$-th collection, $\{\E^i_1,\ldots,\E^i_{n_i}\}_{n_i\geq 1}$, 
form the common boundary of the immersed hypersurfaces $\s_i=\phi_i^\vee$ and 
$\s_{i+1}=\phi_{i+1}^\vee$. 

The hyperplane containing the ellipsoid $\E^i_j$, denoted $H_{P^i_j}$, is tangent to $\Sigma_{\Qbf}$ along $\E^i_j$ and it is the polar of a singular point $P^i_j$ in $\mathcal{F}_{\Qbf}$ {\rm (mentioned in Theorem\,\ref{th:Unavoidable-singularities-FQ}-$c$)}. 
%
\end{theorem}

So these collections of ellipsoids (which consists of parabolic points) separate 
the front $\Sigma_{\Qbf}$ into the 
$\lfloor\frac{k+1}{2}\rfloor$ immersed hypersurfaces 
$\s_i=\phi_i^\vee$ described in Theorem\,\ref{th:configuration-SigmaQ} (with $\s_{\lfloor\frac{k+1}{2}\rfloor}=\eta^\vee$ if $k$ is odd). 

\begin{proof}
An isolated hyperbolic singular point $D$ locally separates $\mathcal{F}_{\Qbf}$ into two open local branches. 
Of course, the duals of these two open branches of $\mathcal{F}_{\Qbf}$ are part of $\Sigma_{\Qbf}$. 
The sets of tangent hyperplanes to these two branches are separated by the $(\ell-2)$-parameter 
family of hyperplanes tangent to $\mathcal{F}_{\Qbf}$ at $D$. These hyperplanes are the tangent hyperplanes 
to the quadratic cone determined by the singular point $D$. By Example\,\ref{ex:ellipse-cone}, 
their poles form an $(\ell-2)$-dimensional ellipsoid $\E_D$ in $\Sigma_{\Qbf}$. Of course, $D$ is the pole of the hyperplane that contains the ellipsoid $\E_D$. 
\end{proof}

\begin{example*}
The four singular points of the Cayley Cubic Surface $\mathcal{F}_{\Qbf}$, in Fig.\,\ref{fig:Cayley-Cubic-Surface}, 
are the poles of the four planes which are tangent to the Steiner Roman Surface $\Sigma_{\Qbf}$, 
in Fig.\,\ref{fig:Steiner-Roman-Surface}, along four respective ellipses 
of parabolic points.   These four ellipses separate $\Sigma_{\Qbf}$ into two immersed 
surfaces $\s_1$ of positive Gauss curvature and $\s_2$ of negative Gauss curvature. 
\end{example*}

\subsection*{B \ $A_1$-singularities in $\mathcal{F}_{\Qbf}$ coorespond to tangent discs in $\Sigma_{\Qbf}$ (if $k>\ell$)} 

\begin{theorem}\label{th:front-singularities}
Let $\Qbf:\R^k\to\R^\ell$ be a homogeneous quadratic map with $k>\ell$. 
\smallskip

\noindent
$(a)$ The front $\Sigma_{\Qbf}$ has no parabolic point if the focal $\mathcal{F}_{\Qbf}$ is smooth. 
\smallskip

\noindent
$(b)$ It the front $\Sigma_{\Qbf}$ is in general position, it has stable  $($Legendrian$)$ front singularities.  
\smallskip

\noindent 
$(c)$ If the focal $\mathcal{F}_{\Qbf}$ has Morse singular points $P_1,\ldots,P_N$, then:
\smallskip

$(c_1)$ The parabolic set of $\Sigma_{\Qbf}$ is the union of $N$ hyperplane closed discs $\Delta_1,\ldots,\Delta_N$,  
$($with multiplicity $2$$)$, consisting of infinitely degenerate parabolic points and bounded by ellipsoids $\E_1,\ldots,\E_N$ of dimension $\ell-2$. 
The hyperplane containing the double disc $\Delta_i$ is the polar of $P_i$ and is tangent 
to $\Sigma_{\Qbf}$ along the ellipsoid $\E_i$. 
\smallskip

$(c_2)$ Under a generic $1$-parameter
deformation of the map $\Qbf$, a neighbourhood of the double disc $2\times\Delta_i$ in $\Sigma_{\Qbf}$ splits into
two hypersurfaces close to the disc $\Delta_i$ - {\rm Fig.\,\ref{fig:ellipsoid-perestroika}}. 
\end{theorem}

\begin{proof}
$(a)$ The parabolic points of $\Sigma_{\Qbf}$ would provide Legendre front singularities in the focal, 
but $\mathcal{F}_{\Qbf}$ does not admit such singularities, due to its hyperbolicity 
(by Th.\,\ref{th:hyperbolic-singularities}). 
\smallskip

\noindent
$(b)$ The Legendre front singularities of $\Sigma_{\Qbf}$ 
correspond by polar duality to the parabolic set of the focal $\mathcal{F}_{\Qbf}$, 
which is stable if $\Qbf$ is in general position. 
\smallskip

\noindent
$(c_1)$ To prove that the doble hyperplane disc $\Delta_i$ bounded by the ellipsoid $\E_i$ is part of $\Sigma_{\Qbf}$, we shall take 
a generic admissible local deformation $\mathcal{F}_{\Qbf}(\e)$ of $\mathcal{F}_{\Qbf}$. 
Near its hyperbolic singular point $P_i$, the deformation  $\mathcal{F}_{\Qbf}(\e)$ has two local components.  

Take one small open disc in each local component of   $\mathcal{F}_{\Qbf}(\e)$. 
Near the singular point $P_i$ 
the respective sets of tangent hyperplanes to these two discs are bounded by non small solid angles (close to the tangent quadratic cone) - Fig.\,\ref{fig:double-pensil-disc}-right. 
If $\mathcal{F}_{\Qbf}(\e)$ tends to $\mathcal{F}_{\Qbf}$, then these bounding solid angles tend to the  
quadratic tangent cone determined by the singular point $P_i$ in $\mathcal{F}_{\Qbf}$ - Fig.\,\ref{fig:double-pensil-disc}-left.  
Thus, in the limit, both sets of tangent hyperplanes become a disc in the pencil of hyperplanes through $P_i$ -  
its boundary is the set of  
hyperplanes tangent to the tangent quadratic cone at $P_i$, whose poles form an $(\ell-2)$-dimensional ellipsoid (Example\,\ref{ex:ellipse-cone}). 
\begin{figure}[ht!]
\centering
\includegraphics [scale=0.4]{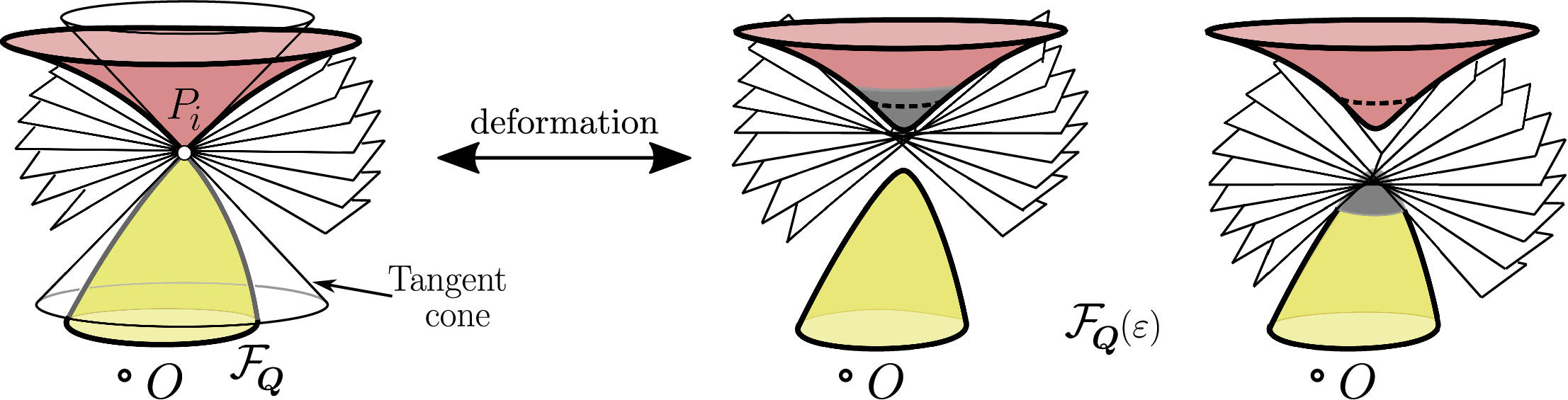}
\caption{\small For each local component of 
$\mathcal{F}_{\Qbf}(\e)$ the set of tangent hyperplanes to a small disc (in gray) tends to a disc of 
the pensil of hyperplanes through the singular point $P_i$ in $\mathcal{F}_{\Qbf}$.}
	\label{fig:double-pensil-disc}
\end{figure}

Hence, besides the duals of the two open local components of $\mathcal{F}_{\Qbf}$, the front 
$\Sigma_{\Qbf}$ contains an ellipsoidal double disc of parabolic points lying in the polar hyperplane of $P_i$. 
Its interior points are the poles of the hyperplanes through $P_i$
that cut the tangent cone only at $P_i$, and its boundary points are the poles of the hyperplanes 
that are tangent to the cone - forming the bounding ellipsoid $\E_i$ - fig.\,\ref{fig:ellipsoid-perestroika}-left.  
\begin{figure}[ht!]
\centering
\includegraphics [scale=0.36]{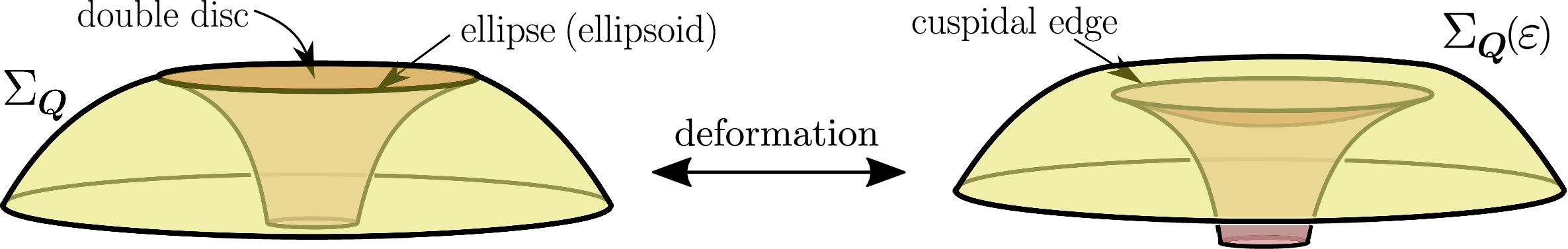}
\caption{\small The front $\Sigma_{\Qbf}$ contains a double ellipsoidal disc $\Delta_i$ 
associated by duality to a singular point $P_i$ in $\mathcal{F}_{\Qbf}$. 
The cuspidal edge in $\Sigma_{\Qbf}(\e)$ is a deformation od the ellipsoid $\E_i$ and is dual to the parabolic set in $\mathcal{F}_{\Qbf}(\e)$} 
\label{fig:ellipsoid-perestroika}
\end{figure}
\smallskip

\noindent
$(c_2)$ A generic deformation $\Qbf(\varepsilon)$ of $\Qbf$ smooths $\mathcal{F}_{\Qbf}$, near 
its singular point $P_i$, into two small local topological discs in $\mathcal{F}_{\Qbf}(\varepsilon)$ - Fig.\,\ref{fig:double-pensil-disc}; then the above double disc of the pencil of hyperplanes 
through $P_i$ splits into two sets of hyperplanes which are tangent to the respective small discs in $\mathcal{F}_{\Qbf}(\varepsilon)$ (Fig.\,\ref{fig:double-pensil-disc}). Then, in the dual $\Sigma_{\Qbf}$, a neighbourhood of the
double disc $2\times\Delta_i$ splits into two hypersurfaces close to that disc - Fig.\,\ref{fig:ellipsoid-perestroika}.
\end{proof}

\begin{figure}[ht!]
\centering
\includegraphics [scale=0.22]{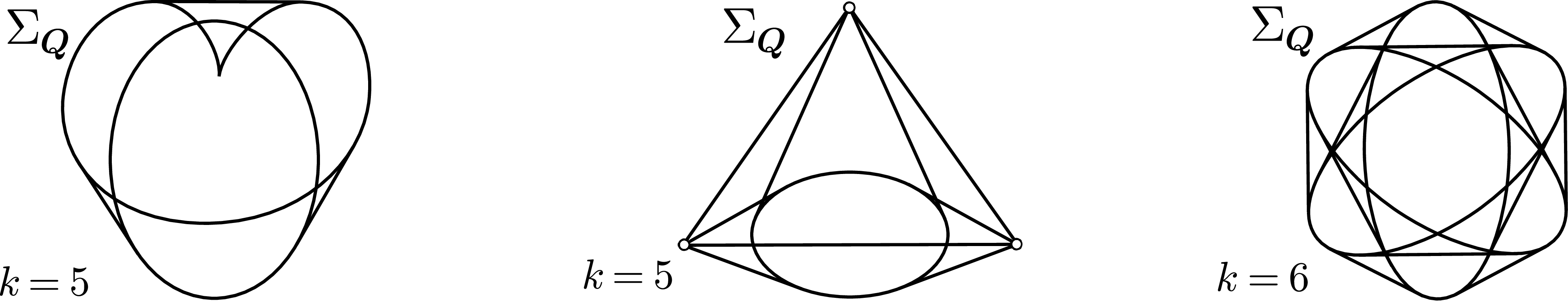}
\caption{\small Fronts 
$\Sigma_{\Qbf}$ of quadratic maps $\Qbf:\R^k\to\R^2$ for $k=5, 5, 6$, duals to the focals of Fig.\ref{fig:nested-top-hyperspheres}}.
\label{fig:fronts_k=5-5-6}
\end{figure}

\begin{example*}
Each double point of a focal $\mathcal{F}_{\Qbf}$ in Fig.\,\ref{fig:nested-top-hyperspheres} 
corresponds to a double bitangent segment of the respective front $\Sigma_{\Qbf}$ in Fig.\,\ref{fig:fronts_k=5-5-6}. 
In Fig.\,\ref{fig:special-quintic}, the double point in the quintic $\mathcal{F}_{\Qbf}$ of Example\,\ref{ex:special-quintic}
corresponds (by duality) to the double segment bitangent to the cardioid in $\Sigma_{\Qbf}$.
\end{example*}

\begin{example*}[\textbf{\textit{Ellipsoidal cuspidal edges in $\Sigma_{\Qbf}$}}]
The hyperbolicity of $\mathcal{F}_{\Qbf}$ implies that 
the quadratic part of its equation has the form $1/2(a_1^2q_1^2+\cdots+a_{\ell-1}^2q_{\ell-1}^2-(q_{\ell}-r)^2)$ with 
$a_i>0$, near an $A_1$-singular point. Take the following deformation $S_\varepsilon$ of an $A_1$-singular point 
$1/2(a_1^2q_1^2+\cdots+a_{\ell-1}^2q_{\ell-1}^2-(q_{\ell}-r)^2)+
\a(q_{\ell}-r)^3/6=-\varepsilon^2$. The parabolic set of 
$S_\e$ is given by the intersection of $S_\e$ with the hypersurface $P$ given by the following equation (obtained from the bordered Hessian): 
$(1-\a(q_{\ell}-r))
(a_1^2q_1^2+\cdots+a_{\ell-1}^2q_{\ell-1}^2-(q_{\ell}-r)^2)
+\a^2(q_{\ell}-r)^4/4=0$. 

Replacing the degree $2$ terms of $P$ in the equation of $S_\e$ we get that the parabolic 
set of $S_\e$ is close to the small ellipsoid $\E$, in the hyperplane of equation $q_{\ell}-r=(-6\e^2/\a)^{1/3}$, 
given by the equation 
$a_1^2q_1^2+\cdots+a_{\ell-1}^2q_{\ell-1}^2=(-6\e^2/\a)^{2/3}$. 

Then the cuspidal edge of $\Sigma_{\Qbf(\e)}$ (see Fig.\,\ref{fig:ellipsoid-perestroika}) is close to the $(\ell-2)$-ellipsoid which consists of the poles tangent to $S_\e$ along $\E$,  (Example\,\ref{ex:ellipse-cone}): 
\[\frac{q_1^2}{a_1^2}+\cdots+\frac{q_{\ell-1}^2}{a_{\ell-1}^2}=\frac{1}{r^2}, \qquad q_\ell=\frac{1}{r}\,.\]

\end{example*}

\section{Homogeneous Quadratic Maps from $\R^k$ to $\R^2$}\label{sect:QuadMapToR^2}  

\subsubsection*{A \ On the Number of Cusps on a Front $\Sigma_{\Qbf}$ in $\R^2$}

\begin{theorem}\label{prop:klein}
Let $\Qbf:\R^k\to\R^2$ a homogeneous quadratic map with $k-1\geq 2$. 
\smallskip

\noindent 
$(a)$ The front $\Sigma_{\Qbf}$ {\rm (critical value set of $\Qp$)} has at most $k(k-2)$ cusps.  
\smallskip

\noindent
$(b)$ There exist quadratic maps $\Qbf:\R^k\to\R^2$ whose front has $k(k-2)$ cusps. 
\end{theorem}

\begin{proof}
$(a)$ On one hand, the cusps of $\Sigma_{\Qbf}$ correspond to the inflections of $\mathcal{F}_{\Qbf}$, which 
is a degree $k$ algebraic curve, and a Pl\"ucker relation states that the number of inflections 
of a non singular algebraic curve of degree $k$ is equal to $3k(k-2)$ (the presence of singularities 
on $\mathcal{F}_{\Qbf}$ would reduce the number of its inflections). 
On the other hand, G.\,Salmon \cite{Salmon-HPC} discovered (and F.\,Klein proved 
\cite{Klein-real-inflections, Klein-elementary3}) that \textit{the number of real inflections of a real algebraic curve is at most one third of the total number of inflections, 
that is, at most $k(k-2)$}. Thus Theorem\,\ref{prop:klein} follows from Salmon-Klein Theorem. 
	
\noindent
$(b)$ The construction of such maps is given in proposition\,\ref{prop:k(k-2)_cusps}, Section\,\ref{ssect:split-DC}\,E. 
\end{proof}

\begin{theorem}\label{min-cusps}
$(a)$ For each even natural number $k=2n$ there are homogeneous quadratic maps $\Qbf:\R^{2n}\to\R^2$ whose front $\Sigma_{\Qbf}$ has no cusp. 
\smallskip

\noindent
$(b)$ If $k$ is odd and $\Qbf:\R^k\to\R^2$ is a homogeneous quadratic map whose focal $\mathcal{F}_{\Qbf}$ is a smooth non reducible curve, then the front $\Sigma_{\Qbf}$ has at least three cusps. 
\smallskip

\noindent
$(c)$ For $k=2n+1$ there are homogeneous quadratic maps $\Qbf:\R^{2n+1}\to\R^2$ whose focal $\mathcal{F}_{\Qbf}$ is a smooth non reducible curve, and whose front $\Sigma_{\Qbf}$ has exactly three cusps.
\end{theorem}
\begin{proof}
$(a)$ For the homogeneous quadratic map $\Qbf:\R^{2n}\to\R^2$ given by 
\begin{equation}\label{eq:map-concentric-circles}
\Qbf(t)=\sum_{i=1}^{n}\frac{n-i+1}{5}\left[
	\begin{pmatrix}
	3 \\
	4
	\end{pmatrix}
(t_{2i-1}^2-t_{2i}^2)+
2\begin{pmatrix}
	4 \\
	-3
\end{pmatrix}
t_{2i-1}t_{2i}
\right]
\end{equation}
the focal $\mathcal{F}_{\Qbf}$ has equation 
\begin{equation}
	\prod_{i=1}^n\left[ 1-\frac{(n-i+1)^2}{5^2}
	\left(\left\langle 
	\begin{pmatrix}
	3 \\
	4
	\end{pmatrix}\,,\,
	\begin{pmatrix}
	q_1 \\
	q_2
	\end{pmatrix}\right\rangle^2  + 
\left\langle 
	\begin{pmatrix}
	4 \\
	-3
	\end{pmatrix}\,,\,
	\begin{pmatrix}
	q_1 \\
	q_2
	\end{pmatrix}\right\rangle^2\,\right)
\right]=0\,,
\end{equation}
which can be written as 
\[
\prod_{i=1}^n\left(q_1^2+q_2^2-\frac{1}{(n-i+1)^2}\right)=0\,.
\]
Thus, $\mathcal{F}_{\Qbf}$ has no inflection - it consists of $n$ circles, centred at the origin, of consecutive radii 
$1/n$, $1/(n-1)$, \ldots, $1$. Therefore $\Sigma_{\Qbf}$ has no cusp and consist of $n$ circles, 
centred at the origin, of decreasing radii $n$, $n-1$, \ldots, $1$.  
\smallskip
	
\noindent
$(b)$ Let us prove that the focal $\mathcal{F}_{\Qbf}$ has at least three inflections.
Theorem\,\ref{th:topological-configurations} states that for odd $k$ one component of the curve 
$\mathcal{F}_{\Qbf}$ is a topological projective line. 
A M\"obius theorem states that \textit{an embedded curve of $\RP^2$ which is a topological projective line 
has at least three inflections}. Thus, by duality, $\Sigma_{\Qbf}$ has at least three cusps. 
\smallskip
	
\noindent
$(c)$ Take the quadratic map $\Qbf:\R^{2n+1}\to\R^2$, given by {\small
		\begin{equation*}\label{eq:map-concentric-circles}
			\Qbf(t)=\sum_{i=1}^{n}\textstyle{\frac{n-i+1}{5}}\left[
			\textstyle{\begin{pmatrix}
					3 \\
					4
			\end{pmatrix}}
			(t_{2i-1}^2-t_{2i}^2)+
			2\textstyle{\begin{pmatrix}
					4 \\
					-3
			\end{pmatrix}}
			t_{2i-1}t_{2i}
			\right]+\textstyle{\frac{2}{5}}
			\textstyle{\begin{pmatrix}
					0 \\
					1
			\end{pmatrix}}
			t_{2n}t_{2n+1}+\textstyle{\frac{4}{5}}
			\textstyle{\begin{pmatrix}
					0 \\
					1
			\end{pmatrix}}
			t_{2n+1}^2\,.
	\end{equation*} }
An easy short calculation leads to the equation of the focal $\mathcal{F}_{\Qbf}$, 
\[
\left[\prod_{i=1}^{n-1}\left(q_1^2+q_2^2-\frac{1}{(n-i+1)^2}\right)\right]
\left[\left(1-q_1^2-q_2^2\right)\left(\textstyle{\frac{4}{5}}q_2-1\right)+
\frac{q_2^2}{5^2}\left( 1-\textstyle{\frac{1}{5}}\left(3q_1+4q_2 \right)\right)
\right]=0\,.
\]
It is a reducible curve which consists of $n-1$ circles centred at the origin (with consecutive increansing radii 
$1/n$, $1/(n-1)$, \ldots , $1/2$) and a smooth cubic whose oval is close to the unit circle and whose second component has three inflections.  
	
Any sufficiently small generic perturbation $\widetilde{\Qbf}$ of $\Qbf$ will have a focal $\mathcal{F}_{\widetilde{\Qbf}}$ close to $\mathcal{F}_{\Qbf}$, non reducible and with no new inflection. Therefore, the front $\Sigma_{\widetilde{\Qbf}}$ consists of $n$ nested ovals and a small curve that has exactly three cusps. 
\end{proof}

\begin{example}[of Theorem\,\ref{min-cusps}\,$a$]
The front $\Sigma_{\Qbf}$ of the quadratic map $\Qbf:\R^4\to\R^2$, with 
\[ \Qbf(t)=\frac{1}{5}
\begin{pmatrix}
	6\left(t_1^2-t_2^2\right)+
	16\,t_1t_2+3\left(t_3^2-t_4^2\right)+
	8\,t_3t_4 \\
	8\left(t_1^2-t_2^2\right)+
	12\,t_1t_2+4\left(t_3^2-t_4^2\right)+
	6\,t_3t_4
\end{pmatrix}\,,  
\]
consist of two circles of respective radii $2$ and $1$, centred at the origin. 
\end{example}

\subsubsection*{B \ Admisible and Forbidden Metamorphosis of the Focal $\mathcal{F}_{\Qbf}$}
\begin{figure}[ht!]
	\centering
	\includegraphics [scale=0.22]{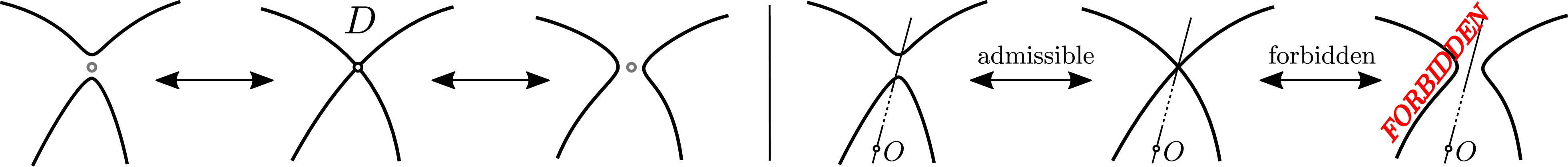}
	\caption{\small Left: Morse bifurcation\,; \hspace*{1.8cm} Right: admissible and forbidden metamorphoses.}
	\label{fig:restricted-bifurcation}
\end{figure}

Since the focal $\mathcal{F}_{\Qbf}$ is not immersed, but implicitly given by an equation, 
its double points may undergo metamorphoses. 
To describe the dual front of a neighbourhood 
of a double point in $\mathcal{F}_{\Qbf}$, 
we shall deform it: 

If a curve given by an equation has a Morse double point, its deformation inside a generic $1$-parameter 
family of curves is a standard Morse bifurcation (Fig.\,\ref{fig:restricted-bifurcation}-left). 

However, the hyperbolicity of a focal $\mathcal{F}_{\Qbf}$ forbids one of the double arrows in 
Fig.\,\ref{fig:restricted-bifurcation}-left (Theorem\,\ref{th:hyperbolicity-theorem}\,a). 
The forbidden double arrow depends on the relative position of the origin. 
Indeed, take the line that joins the origin to the 
double point of $\mathcal{F}_{\Qbf}$. After a local smoothing of $\mathcal{F}_{\Qbf}$, the number of intersection points of every neighbouring line trough the origin with $\mathcal{F}_{\Qbf}$, counting multiplicities, cannot change. So, those lines must cut both local branches. 
Then the left metamorphosis of Fig.\,\ref{fig:restricted-bifurcation}-right is admissible, while the right one is forbidden.   

\subsection*{C \ Double points in $\mathcal{F}_{\Qbf}$ are bitangent (double) segments in $\Sigma_{\Qbf}$} 
\begin{figure}[ht!]
\centering
\includegraphics [scale=0.22]{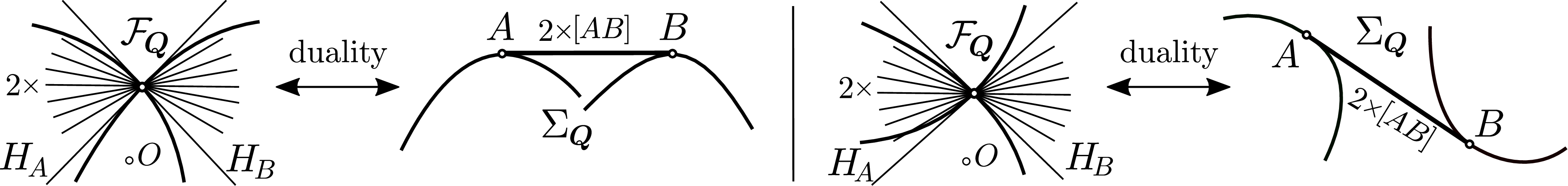}
\caption{\small A segment of pencil at a double point in $\mathcal{F}_{\Qbf}$ is dual  to a bitangent segment in $\Sigma_{\Qbf}$.}
\label{fig:duality-double-point}
\end{figure}

Morse double points in the focal $\mathcal{F}_{\Qbf}$ 
are responsible for the presence of bitangent segments in the front $\Sigma_{\Qbf}$. 
\begin{theorem}\label{th:bitangent-segment}
If the focal $\mathcal{F}_{\Qbf}$ has a Morse double point $D$, then the dual of a neighbourhood 
of $D$ in $\mathcal{F}_{\Qbf}$ is the union of three curves in $\Sigma_{\Qbf}$$:$ the duals of the smooth 
branches of $\mathcal{F}_{\Qbf}$ crossing at $D$, and a segment, with multiplicity two, 
tangent to both dual branches and whose end-points are the respective tangency points 
{\rm - Figs.\,\ref{fig:bitangent-segments} and \ref{fig:duality-double-point}}. 
\begin{figure}[ht!]
\centering
\includegraphics [scale=0.22]{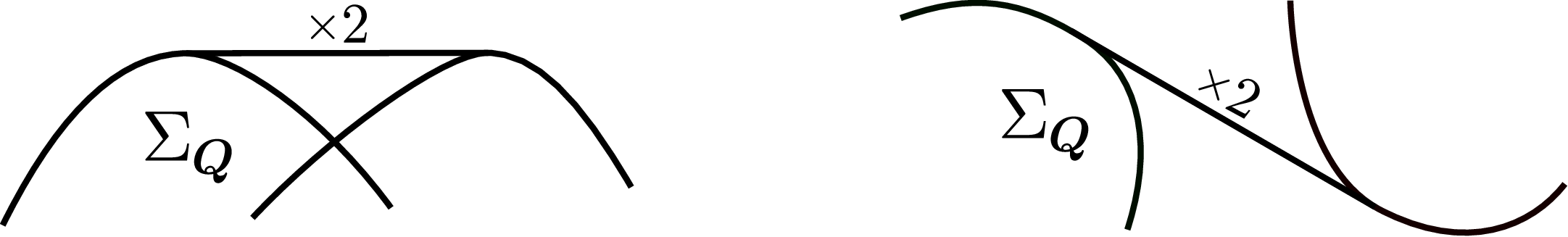}
\caption{\small Possible bitangent segments of the front $\Sigma_{\Qbf}$ (critical value set of $\Qp$).}
\label{fig:bitangent-segments}
\end{figure}
\end{theorem}

\begin{proof}
It is obvious that the duals of the two smooth branches of $\mathcal{F}_{\Qbf}$ crossing at $D$ belong 
to $\Sigma_{\Qbf}$. To prove that the bitangent segment also belongs to $\Sigma_{\Qbf}$, we shall take 
a generic admissible local deformation $\mathcal{F}_{\Qbf}(\e)$ of $\mathcal{F}_{\Qbf}$, near its double 
point $D$ and then describe how $\mathcal{F}_{\Qbf}(\e)$ behaves when it comes 
back to $\mathcal{F}_{\Qbf}$.  

Since the close small arcs of the two local components of $\mathcal{F}_{\Qbf}(\e)$ have hight curvature,   
their respective tangent lines turn quickly along a non small angle - Fig.\,\ref{fig:double-pensil-segment}. 
When $\mathcal{F}_{\Qbf}(\e)$ comes back to $\mathcal{F}_{\Qbf}$ this angle tends to the angle between the smooth 
crossing branches of $\mathcal{F}_{\Qbf}$ and this set of tangent lines becomes a segment of the pencil of lines through 
the double point $D$, taken with multiplicity two (one for each local component) - Fig.\,\ref{fig:double-pensil-segment}. 
\begin{figure}[ht!]
\centering
\includegraphics [scale=0.22]{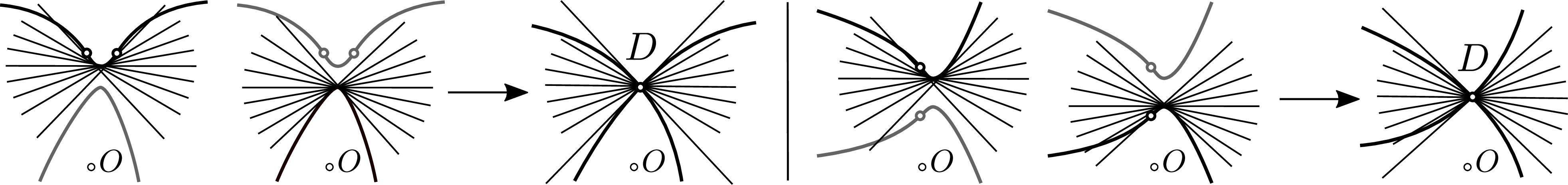}
\caption{\small The neighbouring tangent lines converge to twice a segment of a pencil of lines.}
\label{fig:double-pensil-segment}
\end{figure}

Then the poles of the lines that form that segment of pencil, form a segment $[PQ]$ (taken with multiplicity two) 
of the polar of $D$ whose end-points are the poles of the tangent lines 
to the crossing smooth branches of $\mathcal{F}_{\Qbf}$ - Fig.\,\ref{fig:duality-double-point}. 
\end{proof}


\begin{theorem}\label{th:perestroikas}
Let $\Qbf:\R^k\to\R^2$ a homogeneous quadratic map whose front $\Sigma_{\Qbf}$ contains a 
double bitangent segment $[AB]$ {\rm (see Th.\,\ref{th:bitangent-segment} and Fig.\ref{fig:bitangent-segments})}. 
Under a generic $1$-parameter deformation of $\Qp$, the neighbourhood of the double segment $2\times[AB]$ in $\Sigma_{\Qbf}$ 
splits into two curves close to the segment $[AB]$. There are two possible cases {\rm (Figs.\,\ref{fig:segment-perestroikas} and \ref{fig:duality-perestroika}):} 
\smallskip 

\textit{\textbf{a}}. One curve with two cusps {\rm (near the points $A$ and $B$)} and one smooth curve; 
\medskip

\textit{\textbf{b}}. Both curves have one cusp  {\rm (one cusp near $A$ and the other cusp near $B$).}  
\begin{figure}[ht!]
	\centering
	\includegraphics [scale=0.175]{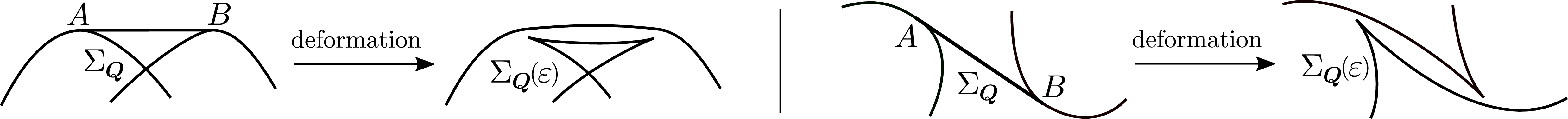}
	\caption{\small The deformation of $\Sigma_{\Qbf}$ near a double bitangent segment (the two possible cases).}
	\label{fig:segment-perestroikas}
\end{figure}
\end{theorem}
\begin{proof}
The points of the double segment $[AB]$ and of its neighbourhood in $\Sigma_{\Qbf}$ are the poles of the tangent lines 
to $\mathcal{F}_{\Qbf}$ in a neighbourhood of a double  point $D$, together with the segment of the pencil of lines 
through $D$ whose ``end-lines'' are the polars of the points $A$ and $B$ - see Fig.\,\ref{fig:duality-double-point}. 
A generic deformation will produce a generic smoothing of the double point of $\mathcal{F}_{\Qbf}$; then  
the double segment of the pencil of lines $2\times[H_AH_B]$ splits into two sets of lines which are tangent to
the respective small arcs of $\mathcal{F}_{\Qbf}$ (Fig.\,\ref{fig:double-pensil-segment}). 
	
At a Morse singular point in general position on a plane curve, six inflection points merge, but only two of them are real 
- the rest being complex (cf. \cite{Plucker, Arnold-V-inflections}). So, in a generic $1$-parameter family, we pass from the case 
where one local branch has two inflections to the case where each local branch has one inflection 
(Fig.\,\ref{fig:restricted-bifurcation}-left). 
	
For the focal $\mathcal{F}_{\Qbf}$, the admissible case depends on the position of 
the origin (see Fig.\,\ref{fig:restricted-bifurcation}-right). 
Thus after a generic deformation $\mathcal{F}_{\Qbf}(\e)$ of a Morse crossing of $\mathcal{F}_{\Qbf}$, the curve $\mathcal{F}_{\Qbf}(\e)$ 
has locally two inflections, and we have two possibilities: 
\smallskip
	
\noindent
a) the two inflections lie on one arc - Fig.\,\ref{fig:duality-perestroika}-left; 
\smallskip
	
\noindent
b) there is one inflection on each arc - Fig.\,\ref{fig:duality-perestroika}-right. 
\begin{figure}[ht!]
	\centering
	\includegraphics [scale=0.215]{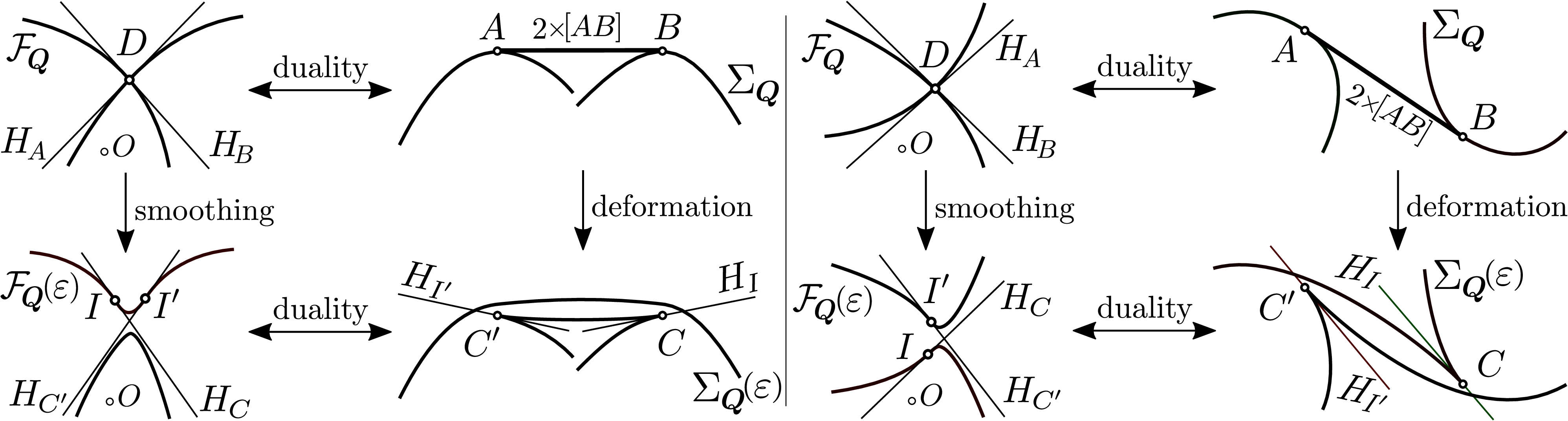}
	\caption{\small The dual of a smoothing transition of a double point in the focal $\mathcal{F}_{\Qbf}$ is the spliting of a double bitangent segment in the front $\Sigma_{\Qbf}$.}
	\label{fig:duality-perestroika}
\end{figure}
	
In both cases, the two inflections of $\mathcal{F}_{\Qbf}(\e)$ are dual to the cusps of $\Sigma_{\Qbf}(\e)$, 
which are respectively close to the end-points $A$, $B$ of the bitangent segment.  
	
Indeed, going back from $\Sigma_{\Qbf}(\e)$ to $\Sigma_{\Qbf}$, the cusps $C, C'$ of $\Sigma_{\Qbf}(\e)$ come to the points $A, B$ at 
which the line $H_D$ is tangent to $\Sigma_{\Qbf}$ because the lines $H_C, H_{C'}$, tangent to $\mathcal{F}_{\Qbf}(\e)$ at its inflections, tend to the lines $H_A, H_B$ tangent to the smooth branches of $\mathcal{F}_{\Qbf}$. 
\end{proof}

\begin{figure}[ht!]
\centering
\includegraphics [scale=0.23]{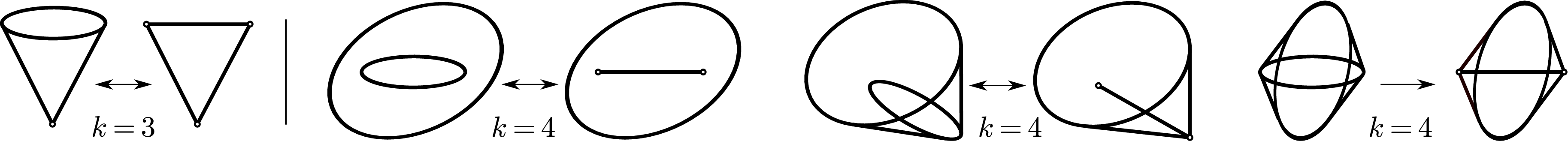}
\caption{\small Exemples of degenerate double segments of $\Sigma_{\Qbf}$ for $\Qbf:\R^k\to\R^2$, $k=3,4$.}
\label{fig:degenerate-segments}
\end{figure}

\begin{remark}[Other Double Segments in the Front]
The front $\Sigma_{\Qbf}$ may have double segments which are not bitangent. 
For example, if the dual $\mathcal{F}_{\Qbf}$ contains two lines intersecting at a double point $D$, 
then the front $\Sigma_{\Qbf}$ contains a segment of the polar line of $D$ (with multiplicity two) -
it is determined by the relative position between the origin and those lines. Under a deformation, 
these lines become a curve close to a hyperbola, while the double segment splits into a curve close to an ellipse - 
Fig.\,\ref{fig:degenerate-segments}. 
\end{remark}

\subsection*{D \ Transitions on the Front $\Sigma_{\Qbf}$ are bitangential events on the Focal $\mathcal{F}_{\Qbf}$}
Besides the spliting of double bitangent segments, 
the fronts $\Sigma_{\Qbf_\e}$ of evolving quadratic maps $\Qbf_\e:\R^k\to\R^2$, 
may undergo the typical transitions of fronts. 
They correspond, via polar duality, to (bi)tangential events on the focal $\mathcal{F}_{\Qbf}$, namely (Fig.\,\ref{fig:front-perestroika2})\,:  

\begin{proposition*}
Let $\Qbf_\e:\R^k\to\R^2$ a generic $1$-parameter family of quadratic maps.  
The front $\Sigma_{\Qbf_\e}$ may ondergo, at isolated parameter values, one of the following transitions$:$ 
\smallskip

\noindent
\textbf{\textit{1}}. Two cusps and a double point are born or die in $\Sigma_{\Qbf}$.
\smallskip

\noindent
\textbf{\textit{2}}. A cusp approaches a smooth branch of $\Sigma_{\Qbf}$, touches it and then crosses that 
branch at two points.
\smallskip

\noindent
\textbf{\textit{3}}. Three smooth branches of $\Sigma_{\Qbf}$ pass through an instant of a triple crossing. 
\begin{figure}[ht!]
\centering
\includegraphics [scale=0.19]{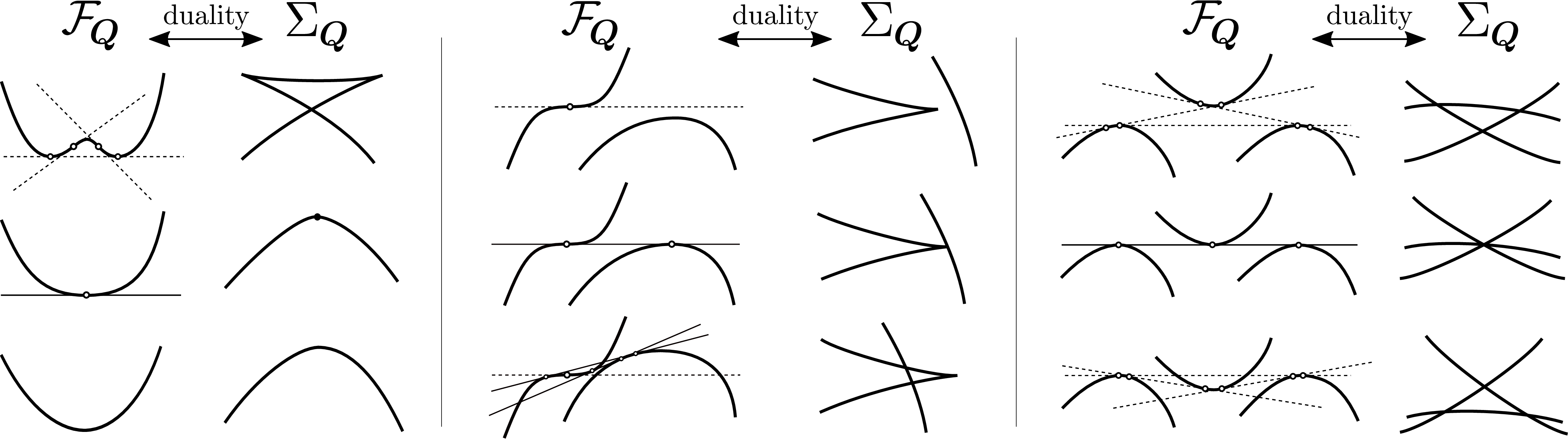}
\caption{\small The three transitions of the front $\Sigma_{\Qbf}$ associated to bitangential events of $\mathcal{F}_{\Qbf}$.}
\label{fig:front-perestroika2}
\end{figure}
\end{proposition*}

\noindent 
Transition 1 
occurs at a point at which $\mathcal{F}_{\Qbf}$ has $4$-point contact with its tangent line; 
Transition 2 occurs if the tangent line to $\mathcal{F}_{\Qbf}$ at an inflection 
is also tangent to another branch of $\mathcal{F}_{\Qbf}$;  
Transition 3 occurs if a line is tangent to three branches of $\mathcal{F}_{\Qbf}$. 

\subsection*{E \ Spliting the Diagonal Case for $\Qbf:\R^k\to\R^2$}\label{ssect:split-DC}
For a diagonal map $\Qbf:\R^k\to\R^2$ the above theorems 
enable us to guess the qualitative shape of a small deformation of the polygon $\Sigma_{\Qbf}$ by a small smoothing 
of the focal curve $\mathcal{F}_{\Qbf}$. Namely, if we smooth $\mathcal{F}_{\Qbf}$ then its double points  
undergo admissible transitions which keep its hyperbolicity 
(Fig.\,\ref{fig:restricted-bifurcation}-right). 

\begin{proposition}
If $\Qbf_\e:\R^k\to\R^2$ is a $1$-parameter family, in general position, of quadratic maps with $\Qbf_0$ a diagonal map, then for $\e$ sufficiently small 
we have$:$
\smallskip 

\noindent
$(a)$ The edges of the polygon $\Sigma_{\Qbf}$ split into two close smooth segments of curve 
in $\Sigma_{\Qbf(\e)}$ {\rm (because all edges of $\Sigma_{\Qbf}$ have multiplicity two)}; 
\smallskip 

\noindent
$(b)$ The deformation $\Sigma_{\Qbf}(\e)$ has no inflection, but may have cusps and crossings; 
\smallskip 

\noindent
$(c)$ All cusps of $\Sigma_{\Qbf}(\e)$ are very close to the vertices of the Polygon $\Sigma_{\Qbf}$ 
{\rm (because the tangents to $\mathcal{F}_{\Qbf}(\e)$ at its inflections are close to the lines that form the curve $\mathcal{F}_{\Qbf}$)};  
\smallskip 

\noindent
$(d)$ The curve $\Sigma_{\Qbf}(\e)$ has at most $k(k-2)$ cusps {\rm (by Theorem\,\ref{prop:klein})}; 
\smallskip 

\noindent
$(e)$ The set $\Sigma_{\Qbf}(\e)$ consists of $k/2$ closed immersed curves if $k$ is even and $(k+1)/2$ closed immersed curves if $k$ is odd. One of these curves is a strictly convex oval that bounds the other closed curves {\rm(by Theorem\,\ref{th:topological-configurations}-$b$)}; 
\smallskip 

\noindent
$(f)$ In a generic bundle of parallel lines there are exactly $k$ distict lines tangent to $\Sigma_{\Qbf}(\e)$ {\rm (by Theorem\,\ref{prop:klein})}. 
\end{proposition}

\begin{proposition}\label{prop:k(k-2)_cusps}
Let $\Qbf_\e:\R^k\to\R^2$ be a small deformation (smoothing) 
of a diagonal map $\Qbf$ whose polygon $\Sigma_{\Qbf}$ is ``convex'' and no three vertices are colinear. Then the front $\Sigma_{\Qbf}(\e)$ has $k(k-2)$ cusps - the maximal possible number.
\end{proposition} 
 
\begin{proof}
By the above six properties of the deformed front $\Sigma_{\Qbf}(\e)$, each of the $k-1$ edges adjacent 
to a given vertex $v$ splits into two segments, and the end-points near $v$ of the resulting segments are joined by pairs 
as follows: the two exterior segments form a smooth arc which is part of the bounding oval (of item $e$) 
and, in order to keep the class of $\Sigma_{\Qbf}(\e)$ equal to $k$ (item $f$) the other $2(k-2)$ segments are joined by pairs into $k-2$ cusps near $v$ 
(Fig.\,\ref{fig:spliting-vertex}). 
Thus each vertex of the polygone $\Sigma_{\Qbf}$ contributes with $k-2$ cusps for the deformed 
front $\Sigma_{\Qbf}(\e)$.  
\begin{figure}[ht!]
\centering
\includegraphics [scale=0.2]{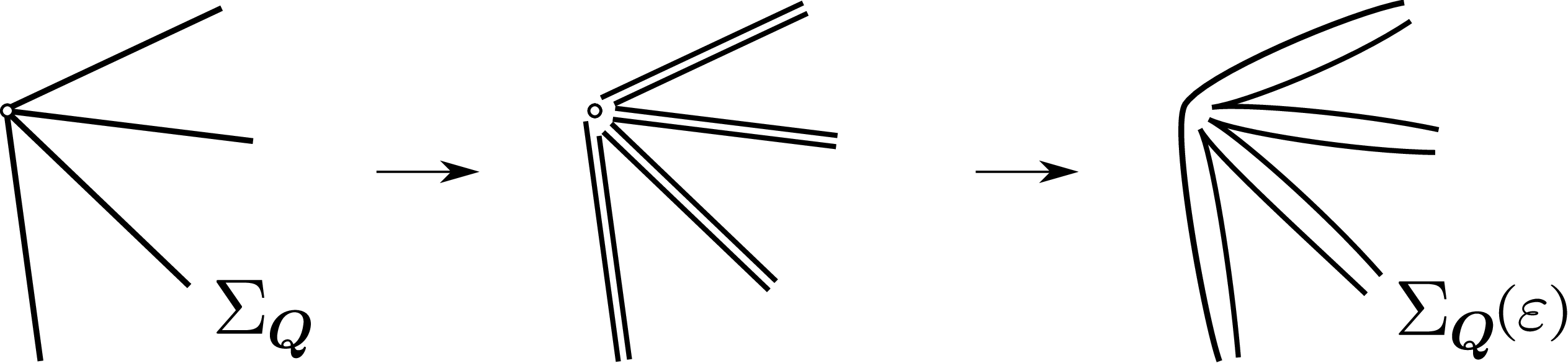}
\caption{\small Spliting the edges adjacent to a vertex of a convex polygon $\Sigma_{\Qbf}$ in the diagonal case.}
\label{fig:spliting-vertex}
\end{figure} 
\end{proof} 

\begin{figure}[ht!]
\centering
\includegraphics [scale=0.157]{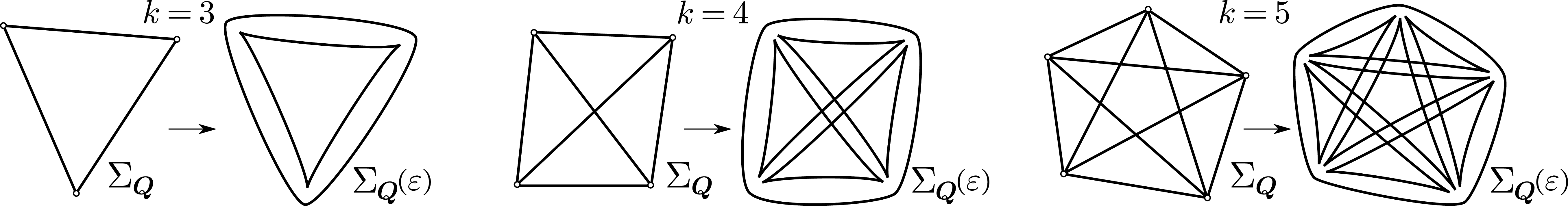}
\caption{\small Spliting a convex polygon $\Sigma_{\Qbf}$ of the diagonal case, by smoothing its focal $\mathcal{F}_{\Qbf}$.}
\label{fig:spliting-diagonal}
\end{figure}

In Fig.\,\ref{fig:spliting-diagonal}, we depict the convex polygons $\Sigma_{\Qbf}$ 
of diagonal maps $\Qbf:\R^k\to\R^2$ and their respective splitings $\Sigma_{\Qbf}(\e)$
(for which $\mathcal{F}_{\Qbf}(\e)$ is smooth) with $k=3,4,5$.  

%
%
%
%
%
%

\section{Classification of Points of $3$-Manifolds in $\R^5$}  

The simplest classification of points of a $k$-manifold $M^k$ in $\R^{k+\ell}$ depend on its order two Taylor expansion which provides the basic local properties of $M^k$. 

Namely, given $p$ in $M$ and writing $\frac{1}{2}\Qbf$ for the local quadratic part (of the Monge form $f:\R^k\to N_pM(=\R^\ell)$) of $M$ at $p$, 
the points are classified in terms of the position of the origin relative to the curvature indicatrix $\Sigma_{\Qbf}$. 

For $M^2$ in $\R^3$ the indicatrix $\Sigma_{\Qbf}$ consists of two points in $N_pM$ (may be coinciding). The are five types of points: elliptic, 
hyperbolic, parabolic, umbilic and flat umbilic.

For $M^2$ in $\R^4$ the local focal $\mathcal{F}_{\Qbf}$ is a conic. The curvature indicatrix $\Sigma_{\Qbf}$ is an ellipse which may degenerate 
into a segment or a point. There are nine possible types of points (only seven of them occur on a surface in general position). 

\begin{figure}[ht!]
\centering
\includegraphics[width=0.98\textwidth]{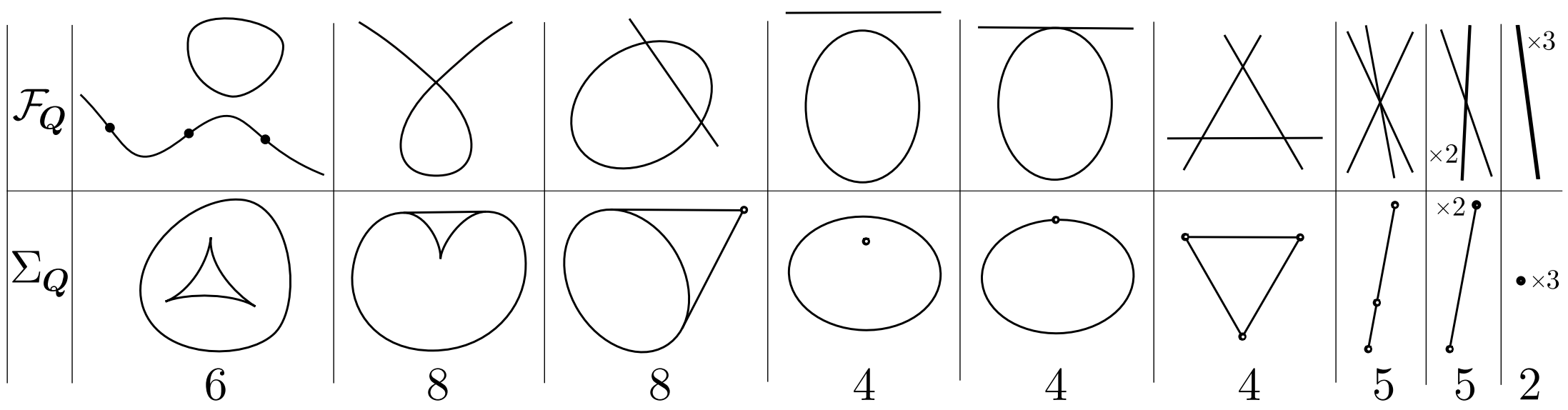} 
\caption{\small Lacal Caustic $\mathcal{F}_{\Qbf}$ and Indicatrix $\Sigma_{\Qbf}$ of a $3$-manifold in $\R^5$. The number of possible types of points is equal to 
$6+8+8+4+4+4+5+5+2=46$.}
\label{fig:cubic-cases}
\end{figure}

\begin{theorem}
For a $3$-manifold $M^3$ in $\R^5$ the local focal $\mathcal{F}_{\Qbf}$ is a cubic curve hyperbolic relative to the origin 
in $\R^2$ - {\rm Fig.\,\ref{fig:cubic-cases}}. The curvature indicatrix $\Sigma_{\Qbf}$ may have the nine possible forms depicted in {\rm Fig.\,\ref{fig:cubic-cases}}. 
There are $46$ possible types of points $($but only $28$ of them occur on a $3$-manifold in $\R^5$ in general position$)$. 
\end{theorem}

{\small 
	\noindent
	\textbf{Acknowledgements}. I'm grateful to Santiago L\'opez de Medrano, for useful discusions and some 
	experimental data. I'm grateful to the Laboratory Solomon Lefschetz UMI2001 CNRS, 
	Universidad Nacional Autonoma de M\'exico, where the main theorem was discovered.}



{\footnotesize 
  \begin{thebibliography}{9}
  


\bibitem {Arnoldwfeeml} {\bf Arnold V.I.}, {\em Wave front evolution and 
equivariant Morse lemma}, Commun. Pure Appl. Math. 29 (1976), 557-582. 

\bibitem {Arnold-V-inflections} {\bf Arnold V.I.}, {\em Vanishing Inflections} 
Functional Analysis and its Appl. 128-130 (1984)

\bibitem {avg} {\bf Arnold V.I., Varchenko A.N., Gussein--Zade S.M.},
{\em Singularities of Differentiable Maps}, Vol. 1, Birkh{\"a}user (1986). 
(Russian version: Nauka, 1982)

\bibitem{AtBG} \textbf{M. F. Atyah M.\,F., Bott\,R., Garding\,L.}, 
\textit{Lacunas for hyperbolic differential operators with constant coefficients. II}, 
Usp.\,Mat.\,Nauk, \textbf{39}, No.\,3 (1984) 171-224.

\bibitem{Garding} \textbf{Garding L.}, \textit{Local hyperbolicity}, 
Isr. J. Math., \textbf{13}, No. 1-2, (1972) 65-81.

\bibitem {Santiago-Samuel} \textbf{Gitler S.}, \textbf{L\'opez de Medrano S.}, 
\textit{Intersections of quadrics, moment-angle manifolds and connected sums}. 
Geometry and Topology \textbf{17} (2013) 1497-1534. 

\bibitem {IRFRT} \textbf{Izumiya S.}, \textbf{Romero Fuster M.C.}, \textbf{Ruas M.A.S.}, 
\textbf{Tari F.}, \textit{Differential Geometry from a Singularity Theory Viewpoint}, 
World Scientific Publishing Co. (2015).

\bibitem {Klein-real-inflections} \textbf{Klein F.}, \textit{Eine neue Relation swischen 
den Singularit\"aten einer algebraischen Kurve}. Mathematische Annalen \textbf{10} (1876) 199-209. 
Reprinted in \textit{Klein's collected works} Vol. II (1922) 78-88. 

\bibitem {Klein-elementary3} \textbf{Klein F.}, 
\textit{Elementary Mathematics from a Higher Standpoint {\rm III:} Precision Mathematics and 
Approximation Mathematics},  Berlin (1928), 
(English translation, Springer 2016).  

\bibitem {Little} {\bf Little J.}, {\em On singularities of submanifolds of higher dimensional 
Euclidean spaces}. Ann. Mat. Pura ed Appl. {\bf 83}:4 (1969) 261-335. 

\bibitem {Plucker} {\bf Pl\"ucker J.}, \textit{Theorie der algebraischen Kurven}, Bonn (1839).

\bibitem {Santiago1} \textbf{L\'opez de Medrano S.}, 
\textit{Singularities of homogeneous quadratic mappings}.  
RACSAM \textbf{108} (2014) 95-112.

\bibitem {Santiago-libro} \textbf{L\'opez de Medrano S.}, 
\textit{Topology and Geometry of Intersections of Ellipsoids in $R^n$}, 
Springer Cham (2023). 

\bibitem {Salmon-HPC} \textbf{Salmon G.}, \textit{A treatise on the Higher Plane Curves},
Hodges and Smith, Dublin (1852). 

\bibitem {Salmon} \textbf{Salmon G.}, \textit{A treatise in analytic geometry of three dimensions},
Chelsea Pub. (1927). 

\bibitem{Ricardo_R2-Rn} \textbf{Uribe-Vargas R.} \textit{On Surfaces of $\R^n$ via Caustics, Gauss Map 
and Pseudo Euclidean Geometry of Quadratic Forms}. Preprint. 

\bibitem{Ricardo_CI-FS-Rn} \textbf{Uribe-Vargas R.} \textit{Differential Geometry of $k$-manifolds in Euclidean Space via Caustics, Polar Duality and Indicatrix Variety} 
Preprint.

\bibitem{Vainshtein} \textbf{Vainshtein A.\,D.},  \textbf{Shapiro B.\,Z.} 
\textit{Singularities on the Boundary of the Hyperbolicity Region}. J. Sov. Math \textbf{52} (1990) 3326-3337. 

\bibitem{Wall} \textbf{Wall C.T.C.} \textit{Stability, pencils and polytopes}. Bull. London
Math. Soc., \textbf{12} (1980), 401­-421.

\end{thebibliography}
}
\end{document}